\documentclass[12pt,fullpage]{article}

\usepackage{amsfonts,graphicx,amsmath,xr,epsfig}

\def\Z{\mathbb{Z}}
\def\R{\mathbb{R}}
\def\N{\mathbb{N}}
\def\epsilon{\varepsilon}

\newcommand{\SE}{\setcounter{equation}{0} \section}
\newcommand{\be}{\begin{equation}}
\newcommand{\ee}{\end{equation}}
\newcommand{\baa}{\begin{array}}
\newcommand{\eaa}{\end{array}}
\newcommand{\bea}{\begin{eqnarray}}
\newcommand{\eea}{\end{eqnarray}}
\newcommand{\beaa}{\begin{eqnarray*}}
\newcommand{\eeaa}{\end{eqnarray*}}

\newtheorem{theo}{\bf Theorem}[section]

\newtheorem{pro}[theo]{\bf Proposition}

\newtheorem{rem}{\bf Remark}[section]

\begin{document}
\date{\today}
\title{\bf{Propagation and blocking in periodically hostile environments}
\thanks{Part of this work was carried out during a visit by the second author at the Department of Mathematics of the National Taiwan Normal University, the hospitality of which is thankfully acknowledged. The first author is partially supported by the National Science Council of the Republic of China under the grant NSC 99-2115-M-032-006-MY3. The second author is supported by the French ``Agence Nationale de la Recherche" within the projects ColonSGS and PREFERED. He is also indebted to the Alexander von Humboldt Foundation for its support.}
}
\author{Jong-Shenq Guo$^{\hbox{\small{ a}}}$ and Fran\c cois Hamel$^{\hbox{\small{ b}}}$
\\
\\
\footnotesize{$^{\hbox{a }}$Tamkang University, Department of Mathematics}\\
\footnotesize{151, Ying-Chuan Road, Tamsui, New Taipei City 25137, Taiwan}\\
\footnotesize{$^{\hbox{b }}$Aix-Marseille Universit\'e \& Institut Universitaire de France}\\
\footnotesize{LATP, Facult\'e des Sciences et Techniques, F-13397 Marseille Cedex 20, France}}
\maketitle

\begin{abstract}
We study the persistence and propagation (or blocking) phenomena for a species
in periodically hostile environments.
The problem is described by a reaction-diffusion equation with zero Dirichlet boundary condition.
We first derive the existence of a minimal nonnegative nontrivial stationary solution and study
the large-time behavior of the solution of the initial boundary value problem.
To the main goal, we then study a sequence of approximated problems in the whole space with
reaction terms which are with very negative growth rates outside the domain under investigation.
Finally, for a given unit vector, by using the information of the minimal speeds of approximated problems,
we provide a simple geometric condition for the blocking of propagation
and we derive the asymptotic behavior of the approximated pulsating travelling fronts.
Moreover, for the case of constant diffusion matrix,
we provide two conditions for which the limit of approximated minimal speeds is positive.
\end{abstract}


\SE{Introduction and main results}\label{intro}

This paper is concerned with persistence and propagation phenomena
for reaction-diffusion equations of the type \be\label{eq}
u_t-\nabla\cdot(A(x)\nabla u)=f(x,u) \ee in $\R^N$ or in some
unbounded open subsets~$\Omega$ of~$\R^N$ with zero Dirichlet
boundary condition on~$\partial\Omega$. Equations of the
type~(\ref{eq}) arise especially in population dynamics and
ecological models (see e.g. \cite{mu,sk,x6}), where the nonnegative
quantity~$u$ typically stands for the concentration of a species.

Let us start with the case of the whole space $\R^N$. The symmetric
matrix field $x\mapsto A(x)=(A_{ij}(x))_{1\le i,j\le N}$ is assumed
to be of class~$C^{1,\alpha}(\R^N)$ with $\alpha>0$ and uniformly
positive definite: that is, there exists a positive
constant~$\beta>0$ such that \be\label{A}
\forall\,x\in\R^N,\quad\forall\,\xi=(\xi_1,\ldots,\xi_N)\in\R^N,\quad
A\xi\cdot\xi:=\sum_{1\le i,j\le
N}A_{ij}(x)\xi_i\xi_j\ge\beta\,|\xi|^2, \ee where $|\cdot|$ denotes
the Euclidean norm in $\R^N$. We set $\R_+=[0,+\infty)$. The
nonlinear reaction term $f:\R^N\times\R_+\to\R,\ (x,u)\mapsto
f(x,u)$ is assumed to be continuous, of class~$C^{0,\alpha}$ with
respect to~$x$ locally uniformly in $u\in\R_+$, of class~$C^1$ with
respect to~$u$, and $\frac{\partial f}{\partial u}(\cdot,0)$ is of
class~$C^{0,\alpha}(\R^N)$. Furthermore, we assume that
\be\label{f}
\left\{\baa{l}
\mbox{$f(x,0)=0$ for all $x\in\R^N$,}\vspace{3pt}\\
\mbox{there exists $M>0$ such that $f(x,M)\le 0$ for all
$x\in\R^N$.}\eaa\right.
\ee
The functions $A_{ij}$ (for all $1\le
i,j\le N$) and $f(\cdot,u)$ (for all $u\in\R_+$) are assumed to be
periodic in $\R^N$. Hereafter a function $w$ is called periodic
in~$\R^N$ if it satisfies
$$w(\cdot+k)=w(\cdot)\ \hbox{ for all }k\in L_1\Z\times\cdots\times L_N\Z,$$
where $L_1,\cdots, L_N$ are some positive real numbers, which are fixed throughout this paper.

If $f$ fulfills the additional Fisher-KPP (for Kolmogorov, Petrovsky and Piskunov)~\cite{f,kpp} assumption
\be\label{fkppn}
\forall\ x\in\R^N,\quad u\mapsto g(x,u)=\frac{f(x,u)}{u}\hbox{ is decreasing with respect to }u>0,
\ee
then the large-time behavior of the solutions of the Cauchy problem
\be\label{cauchyn}\left\{\baa{rcll}
u_t-\nabla\cdot(A(x)\nabla u) & = & f(x,u), & t>0,\ x\in\R^N,\vspace{3pt}\\
u(0,x) & = & u_0(x), & x\in\R^N\eaa\right. \ee is directly related
to the sign of the principal periodic eigenvalue~$\lambda_{1}$ of
the linearized operator at~$0$ (see~\cite{bhr1}). This
eigenvalue~$\lambda_{1}$ is characterized by the existence of a
(unique up to multiplication) periodic function~$\varphi\in
C^{2,\alpha}(\R^N)$, which satisfies
\be\label{eqvarphi}\left\{\baa{rcll}
-\nabla\cdot(A(x)\nabla\varphi)-\zeta(x)\varphi & = & \lambda_{1}\,\varphi & \hbox{in }\R^N,\vspace{3pt}\\
\varphi & > & 0 & \hbox{in }\R^N,\eaa\right. \ee where
$\zeta(x)=\frac{\partial f}{\partial u}(x,0)$ for all $x\in\R^N$.
The precise statement of what is known under the additional
assumption~(\ref{fkppn}) will be recalled just after
Proposition~\ref{pro1} below.

Our first result, which is a preliminary step before the main
purpose of the paper devoted to propagation phenomena in
environments with hostile boundaries, is actually concerned with the
existence of a minimal positive stationary solution~$p$ for
problem~(\ref{cauchyn}) and with the large time behavior of the
solutions $u$ of~(\ref{cauchyn}), when~$f$ fulfills the
assumption~(\ref{f}) alone.

\begin{pro}\label{pro1}
Assume that $\lambda_{1}<0$ and \eqref{f}. Then there is a minimal periodic solution~$p(x)$ of
\be\label{station}\left\{\baa{rcll}
-\nabla\cdot(A(x)\nabla p) & = & f(x,p(x)) & \hbox{in }\R^N,\vspace{3pt}\\
p & > & 0 & \hbox{in }\R^N,\eaa\right. \ee in the sense that, for
any solution $q$ of $(\ref{station})$, there holds $q\ge p$ in
$\R^N$. Furthermore,~$p\le M$ in $\R^N$ and, if $u_0:\R^N\to[0,M]$
is uniformly continuous and not identically $0$, then the solution
$u(t,x)$ of the Cauchy problem~$(\ref{cauchyn})$ is such that
$$\liminf_{t\to+\infty}\,u(t,x)\ge p(x)\ \hbox{ locally uniformly with respect to }x\in\R^N.$$
If one further assumes that $u_0\le p$ in $\R^N$, then $u(t,x)\to
p(x)$ as $t\to+\infty$ locally uniformly with respect to $x\in\R^N$.
\end{pro}

It is obvious to see that the solution $p$ of (\ref{station}) is not
unique in general. Choose for instance $A(x)=I_N$ (the identity
matrix) and $f(x,u)=\sin(u)$ for all $(x,u)\in\R^N\times\R_+$: the
function~$f$ satisfies (\ref{f}) with $M=\pi$, $\lambda_{1}=-1<0$,
but any constant function $p(x)=m\pi$ with $m\in\N\backslash\{0\}$
solves~(\ref{station}). On the other hand, if, in addition
to~(\ref{f}), the function~$f$ satisfies the
assumption~(\ref{fkppn}), then the solution $p$ of (\ref{station})
is unique, see \cite{bhr1}. In particular, all solutions
of~(\ref{station}) are necessarily periodic. Notice that, in the
general case of assumption~(\ref{f}) alone, Proposition~\ref{pro1}
still states the existence of a minimal periodic solution~$p$
of~(\ref{station}) in the class of all positive solutions~$q$, which
are not a priori assumed to be periodic. It is also known that,
under hypotheses~(\ref{f}) and~(\ref{fkppn}), the condition
$\lambda_{1}<0$ of the unstability of $0$ is a necessary condition
for the existence of the solution~$p$ of~(\ref{station}) as well: if
$\lambda_{1}\ge 0$, then all bounded solutions~$u$
of~(\ref{cauchyn}) converge to $0$ as $t\to+\infty$ uniformly in
$\R^N$, see \cite{bhr1}. On the other hand, under the
assumptions~(\ref{f}),~(\ref{fkppn}) and $\lambda_{1}<0$, for any
non-zero bounded uniformly continuous~$u_0:\R^N\to\R_+$, there holds
$u(t,x)\to p(x)$ as $t\to+\infty$ locally uniformly in $x\in\R^N$
(see~\cite{bhr1,fg,w}). We also refer to \cite{bhr1,cc1,cc2} for
related results in the case of bounded domains with Dirichlet or
Neumann boundary conditions, and to \cite{bhna,bhr3} for results
with KPP nonlinearities and periodic or non-periodic coefficients in
$\R^N$. Lastly, it is worth noticing that Proposition~\ref{pro1} and
the aforementioned convergence results are different from what
happens with other types of nonlinearities $f$, like combustion,
bistable or even monostable nonlinearities which are degenerate at
$0$: in these cases, the large-time behavior of the solutions~$u$
of~(\ref{cauchyn}) strongly depends on some threshold parameters
related to the size and/or the amplitude of $u_0$ (see e.g.
\cite{aw,dm,p,r,z1}).

\begin{rem}{\rm The assumption that $u_0$ ranges in the interval $[0,M]$ is made to guarantee
the global existence and boundedness (from below by $0$ and from
above by $M$) of the solutions~$u$ of the Cauchy
problem~(\ref{cauchyn}). If~$f$ fulfills the KPP
assumption~(\ref{fkppn}) together with~(\ref{f}), or if $f(x,s)\le
0$ for all $(x,s)\in\R^N\times[M,+\infty)$, then it follows that the
solution~$u$ exists for all $t\ge 0$ and is globally bounded from
below by $0$ and from above by
$\max\big(M,\|u_0\|_{L^{\infty}(\R^N)}\big)$, as long as $u_0$ is
nonnegative and bounded. The same comment also holds for the Cauchy
problem (\ref{cauchy}) below with zero Dirichlet boundary condition
on $\partial\Omega$.}
\end{rem}

As a matter of fact, in Proposition~\ref{pro1}, the negativity of
$\lambda_{1}$ immediately implies that the positive periodic
functions~$\epsilon\,\varphi$ are subsolutions of~(\ref{cauchyn})
for~$\epsilon>0$ small enough, where~$\varphi$ is a solution
of~(\ref{eqvarphi}). It then follows from the above proposition and
the results of Weinberger~\cite{w} that, for each unit vector~$e$
of~$\R^N$, there is a positive real number $c^*(e)>0$ (minimal
speed) such that the following holds: for each $c\ge c^*(e)$, there
is a pulsating travelling front
$$u(t,x)=\phi(x\cdot e-ct,x)$$
solving $(\ref{cauchyn})$ and connecting $0$ to $p$, that is, the
function $\phi:\R\times\R^N\to[0,M],\ (s,x)\mapsto\phi(s,x)$ is
periodic in $x$, decreasing in~$s$, and it satisfies
$\phi(-\infty,x)=p(x)$ and $\phi(+\infty,x)=0$ for all $x\in\R^N$.
Furthermore, such pulsating travelling fronts do not exist for any
$c<c^*(e)$. We also refer to \cite{n1,nrx,nx,x3,x6} for other
results about pulsating travelling fronts in the whole space~$\R^N$,
including other types of nonlinearities and the case of
time-periodic media.

Now, based on the previous results in $\R^N$, we turn our attention
to the main concern of this paper, namely the case when there are
hostile periodic patches in the domain under consideration. We deal
with persistence and propagation phenomena for reaction-diffusion
equations of the type \be\label{eqper}\left\{\baa{rcll}
u_t-\nabla\cdot(A(x)\nabla u) & = & F(x,u), & x\in\overline{\Omega},\vspace{3pt}\\
u(t,x) & = & 0, & x\in\partial\Omega,\eaa\right. \ee in an unbounded
open set $\Omega\subset\R^N$ which is assumed to be of class
$C^{2,\alpha}$ (with $\alpha>0$) and periodic. The periodicity means
that $\Omega=\Omega+k$  for all $k\in L_1\Z\times\cdots\times
L_N\Z$.  Furthermore, the fields $A(x)$ and~$F(x,u)$ are assumed to
be periodic with respect to $x$ in $\overline{\Omega}$, to have the
same smoothness as before and to fulfill~(\ref{A}) and~(\ref{f})
above, where~$x\in\R^N$ is now replaced
with~$x\in\overline{\Omega}$. In particular, assumption~(\ref{f}) is
now replaced with \be\label{ff}\left\{\baa{l}
\mbox{$F(x,0)=0$ for all $x\in\overline{\Omega}$,}\vspace{3pt}\\
\mbox{there exists $M>0$ such that $F(x,M)\le 0$ for all $x\in\overline{\Omega}$.}\eaa\right.
\ee
Throughout the paper, we denote
$$C=\overline{\Omega}\cap\big(\,[0,L_1]\times\cdots\times[0,L_N]\,\big)$$
the cell of periodicity of $\overline{\Omega}$. The zero Dirichlet
boundary condition imposed on~$\partial\Omega$ means that the
boundary is lethal for the species. Note that the unbounded periodic
open set~$\Omega$ is not a priori assumed to be connected. The
reason for that will become clear later, once the approximation
procedure~(\ref{hypfn}) below has been introduced. However, due to
the global smoothness of~$\partial\Omega$, the set~$\Omega$ has only
a finite number of connected components relatively to the
lattice~$L_1\Z\times\cdots\times L_N\Z$. That is, there is a finite
number of connected components~$\omega_1,\ldots,\omega_m$
of~$\Omega$ such that $\omega_i\cap(\omega_j+k)=\emptyset$ for all
$1\le i\neq j\le m$ and for all $k\in L_1\Z\times\cdots\times
L_N\Z$, and \be\label{Omegai} \Omega=\mathop{\bigcup}_{1\le i\le
m}\Omega_i,\ \ \hbox{ where }\Omega_i=\mathop{\bigcup}_{k\in
L_1\Z\times\cdots\times L_N\Z}\omega_i+k. \ee The sets $\omega_i$
are not uniquely defined, but the sets $\Omega_i$ are unique (up to
permutation), periodic, and $\Omega_i\cap\Omega_j=\emptyset$ for all
$1\le i\neq j\le m$.

In the case of no-flux boundary conditions $\nu(x)\cdot(A(x)\nabla u(t,x))=0$
on~$\partial\Omega$ when $\Omega$ is connected, much work have been devoted in the recent years to the
study of propagation of pulsating fronts $u(t,x)=\phi(x\cdot e-ct,x)$, where~$\phi(s,\cdot)$ is
periodic for all $s\in\R$ and~$e$ is any unit vector, for various types of nonlinearities~$F$,
in straight infinite cylinders~\cite{bn,r} or in periodic domains~\cite{bh,h8,hr,lz1,w}.
In the case of KPP nonlinearities~$F$, further properties of the minimal propagation speeds can be
found in \cite{bhr2,e1,h,kks,llm,n4,rz,z2}.

In this paper, we consider a larger class of reaction terms $F$,
together with zero Dirichlet boundary condition. Let us first
mention that, under the assumption that the equation~(\ref{eqper})
is invariant in the direction~$x_1$ and under appropriate conditions
on~$F$, classical travelling fronts
$$u(t,x)=\phi(x_1-ct,x_2,\ldots,x_n)$$
in straight infinite cylinders (in the $x_1$-direction) with zero
Dirichlet boundary condition are known to exist (see \cite{mn,v2},
including the case of some systems of equations). In this case, the
profiles~$\phi$ of these travelling fronts solve elliptic equations
or systems. For pro\-blem~(\ref{eqper}) in periodic domains, the
reduction to elliptic equations does not hold anymore since the
equation is not assumed to be invariant in any direction. Recently,
existence results for problems of the type~(\ref{eqper}) in
connected two-dimensional periodically oscillating  infinite
cylinders with homogeneous isotropic diffusion ($A(x)=I_2$) and KPP
nonlinearities satisfying~(\ref{fkppn}) have been established,
see~\cite{lz2}. In the present paper, the set~$\Omega$ is periodic
in all variables $x_1,\ldots,x_N$ and the direction of propagation
may be any unit vector~$e$ of~$\R^N$. Actually, one of the novelties
of this paper with respect to the previous literature is that the
nature of propagation vs. blocking strongly depends on the
direction~$e$ and on geometrical properties of the set~$\Omega$
itself.

Let $\lambda_{1,D}$ denote the principal periodic eigenvalue of the
linearized equation at $0$ in $\overline{\Omega}$ with zero
Dirichlet boundary condition. That is, there exists a function
$\varphi\in C^{2,\alpha}(\overline{\Omega})$, which is periodic in
$\overline{\Omega}$ and satisfies
\be\label{lambda1D}\left\{\baa{rcll}
-\nabla\cdot(A(x)\nabla\varphi)-\zeta(x)\varphi & = & \lambda_{1,D}\,\varphi & \hbox{in }\overline{\Omega},\vspace{3pt}\\
\varphi & = & 0 & \hbox{on }\partial\Omega,\vspace{3pt}\\
\varphi & \ge & 0 & \hbox{in }\overline{\Omega},\vspace{3pt}\\
\displaystyle{\mathop{\max}_{\overline{\Omega}}}\,\varphi & > & 0, &
\eaa\right. \ee where $\zeta(x)=\frac{\partial F}{\partial u}(x,0)$
for all $x\in\overline{\Omega}$. If~$\Omega$ is connected,
then~$\varphi>0$ in~$\Omega$ and~$\varphi$ is unique up to
multiplication. Otherwise, in the general case, the
function~$\varphi$ is unique up to multiplication in each set
$\Omega_i$ on which it is positive. More precisely,~$\varphi$ can be
chosen to be positive on the (largest possible) set
$\widetilde{\Omega}=\bigcup_{i\in I_{min}}\Omega_i$, where $I_{min}$
denotes the set of indices~$i\in\{1,\ldots,m\}$ for which the
principal periodic eigenvalue~$\lambda_{1,\Omega_i,D}$ of the
operator $-\nabla\cdot(A(x)\nabla)-\zeta(x)$ in~$\Omega_i$ with zero
Dirichlet boundary condition on~$\partial\Omega_i$ is equal
to~$\lambda_{1,D}$. That is,
$$\lambda_{1,D}=\min_{1\le j\le m}\lambda_{1,\Omega_j,D}=\lambda_{1,\Omega_i,D}\hbox{ for all }i\in I_{min}.$$

The following theorem, which is analogue to Proposition~\ref{pro1},
is concerned with the existence of a minimal nonnegative and
non-trivial stationary solution of~(\ref{eqper})
in~$\overline{\Omega}$ and the large-time behavior of the solutions
of the associated initial boundary value problem, under the
assumption that the steady state~$0$ of~(\ref{eqper}) is linearly
strictly unstable. To do so, we introduce the set \be\label{I-}
I_-=\Big\{i\in\{1,\ldots,m\},\ \ \lambda_{1,\Omega_i,D}<0\Big\}. \ee

\begin{theo}\label{th2}
Assume that $\lambda_{1,D}<0$, that is $I_-\neq\emptyset$. Then
there exists a minimal stationary periodic solution $p(x)$ of
\be\label{statio}\left\{\baa{rcll}
-\nabla\cdot(A(x)\nabla p) & = & F(x,p(x)) & \hbox{in }\overline{\Omega},\vspace{3pt}\\
p & = & 0 & \hbox{in }\partial\Omega\,\cup\,\bigcup_{i\not\in I_-}\!\!\Omega_i,\vspace{3pt}\\
p & > & 0 & \hbox{in }\bigcup_{i\in I_-}\!\!\Omega_i,\eaa\right.
\ee
in the sense that any bounded solution $q$ of $(\ref{statio})$ satisfies $q\ge p$ in $\overline{\Omega}$.
Moreover, for any uniformly continuous function $u_0:\overline{\Omega}\to[0,M]$ which is not identically $0$,
the solution $u(t,x)$ of the initial boundary value problem
\be\label{cauchy}\left\{\baa{rcll}
u_t-\nabla\cdot(A(x)\nabla u) & = & F(x,u), & t>0,\ x\in\overline{\Omega},\vspace{3pt}\\
u(t,x) & = & 0, & t>0,\ x\in\partial\Omega,\vspace{3pt}\\
u(0,x) & = & u_0(x), & x\in\Omega\eaa\right. \ee is such that
\be\label{up} \liminf_{t\to+\infty}\,u(t,x)\ge p(x) \ee locally
uniformly with respect to the points $x\in\overline{\Omega}$ whose
connected components intersect the support of $u_0$. If one further
assumes that $u_0\le p$ in $\overline{\Omega}$, then $u(t,x)\to
p(x)$ as $t\to+\infty$ in the same sense as above.
\end{theo}

As already emphasized, the periodic open set $\Omega$ is not assumed
to be connected, this is why the lower bound (\ref{up}) or the
convergence of $u(t,x)$ to $p(x)$ at large time can only hold in the
(open) connected components~$\mathcal{C}$ of the intersection of
$\Omega$ with the support of $u_0$ (outside these components, the
solution $u(t,x)$ stays $0$ for all times $t\ge 0$). If such a
connected component~$\mathcal{C}$ is included in a set~$\Omega_i$
with $i\in I_-$, then Theorem~\ref{th2} implies that~$u(t,x)$ is
separated away from $0$ at large time, locally uniformly
in~$\mathcal{C}$. However,~(\ref{up}) does not say anything about
the behavior of~$u(t,x)$ when $x\in\bigcup_{i\not\in
I_-}\overline{\Omega_i}$ ($p(x)=0$ there). Actually, for each
$\Omega_i$ with $i\not\in I_-$, one has $\lambda_{1,\Omega_i,D}\ge
0$ and if $F$ satisfies the additional assumption~(\ref{fkppn}) in
$\Omega_i$, then $u(t,x)\to 0$ as $t\to+\infty$ uniformly in
$x\in\overline{\Omega_i}$, as follows from the same ideas as
in~\cite{bhr1}.

The remaining part of this paper is concerned with the existence of
pulsating fronts and the possibility of blocking phenomena for
problem (\ref{eqper}) with zero Dirichlet boundary condition. The
strategy, which is one of the main interests of the paper, consists
in approximating the Dirichlet condition on $\partial\Omega$ (and
even in~$\R^N\backslash\Omega$) by reaction terms with very negative
growth rates in~$\R^N\backslash\overline{\Omega}$, using the
previous results and then passing to the singular limit in the
stationary solutions and in the pulsating travelling fronts as the
growth rates converge to~$-\infty$
in~$\R^N\backslash\overline{\Omega}$. This means that the
quantity~$u$ lives in the whole space~$\R^N$, but the space contains
very bad regions. We will see that the location of the good vs. bad
regions plays a crucial role in the dynamical behavior of the
solutions.

For this, let $(f_n)_{n\in\N}$ be a sequence of real-valued
functions defined in $\R^N\times\R_+$ such that each function
$f_n:(x,u)\mapsto f_n(x,u)$ is continuous, periodic with respect to
$x\in\R^N$, of class~$C^{0,\alpha}$ with respect to $x\in\R^N$
locally uniformly in $u\in\R_+$, of class $C^1$ with respect to~$u$
with $\frac{\partial f_n}{\partial u}(\cdot,0)\in
C^{0,\alpha}(\R^N)$, and it satisfies (\ref{f}). Here we define $\N$ to be the set of all nonnegative integers. Furthermore, we
assume that \be\label{hypfn}\left\{\baa{ll}
f_n(x,u)=F(x,u) & \hbox{ for all }(x,u)\in\overline{\Omega}\times\R_+\hbox{ and }n\in\N,\vspace{3pt}\\
(f_n(x,u))_{n\in\N}\hbox{ is nonincreasing} & \hbox{ for all }(x,u)\in\overline{\Omega}\times\R_+,\vspace{3pt}\\
g_n(x,u)\to-\infty\hbox{ as }n\to+\infty & \hbox{ locally uniformly in }
(x,u)\in(\R^N\backslash\overline{\Omega})\times\R_+,\eaa\right.
\ee
where
$$g_n(x,u)=\left\{\baa{ll}\displaystyle{\frac{f_n(x,u)}{u}} & \hbox{if }u>0,\vspace{3pt}\\
\displaystyle{\frac{\partial f_n}{\partial u}}(x,0)=:\zeta_n(x) &
\hbox{if }u=0.\eaa\right.$$ The last condition means that the death
rate in the region $\R^N\backslash\overline{\Omega}$ is very high,
namely this region becomes more and more unfavorable for the species
as~$n$ becomes larger and larger.

Typical examples of such functions $f_n$ satisfying (\ref{ff}) and (\ref{hypfn}) are
$$f_n(x,u)=\rho_n(x)\,u+\widetilde{f}(u),$$
where the function $\widetilde{f}:\R_+\to\R$ is of class $C^1$ and
satisfies $\widetilde{f}(0)=0$, $\widetilde{f}(M)\le 0$, and the
functions $\rho_n:\R^N\to\R$ are periodic, nonpositive, of class
$C^{0,\alpha}(\R^N)$, nonincreasing with respect to $n$, independent
of $n$ in $\Omega$, and $\rho_n\to-\infty$ as $n\to+\infty$ locally
uniformly in $\R^N\backslash\overline{\Omega}$.

For every $n\in\N$, let $\lambda_{1,n}$ denote the principal
periodic eigenvalue of the linearized operator at~$0$ in $\R^N$.
That is, there exists a (unique up to multiplication) periodic
function~$\varphi_n$ of class~$C^{2,\alpha}(\R^N)$, which satisfies
\be\label{lambda1n}\left\{\baa{rcll}
-\nabla\cdot(A(x)\nabla\varphi_n)-\zeta_n(x)\varphi_n & = & \lambda_{1,n}\,\varphi_n & \hbox{in }\R^N,\vspace{3pt}\\
\varphi_n & > & 0 & \hbox{in }\R^N.\eaa\right. \ee We first
establish the relationship between the principal
eigenvalues~$\lambda_{1,n}$ of~(\ref{lambda1n}) and the principal
eigenvalue~$\lambda_{1,D}$ of~(\ref{lambda1D}), as well as the
convergence of the minimal solutions~$p_n$ of~(\ref{station}) with
nonlinearities~$f_n$ to the minimal solution~$p$ of~(\ref{statio}),
when~$\lambda_{1,D}<0$.

\begin{theo}\label{th3}
Under the above notation, the sequence $(\lambda_{1,n})_{n\in\N}$
is nondecreasing and there holds $\lambda_{1,n}\to\lambda_{1,D}$ as
$n\to+\infty$. Furthermore, if~$\lambda_{1,D}<0$, then the
sequence~$(p_n)_{n\in\N}$ of minimal solutions of~$(\ref{station})$
with nonlinearities~$f_n$ is nonincreasing and
$$p_n(x)\to p_{\infty}(x)\hbox{ as }n\to+\infty\hbox{ for all }x\in\R^N,$$
where, up to a negligible set, $p_{\infty}$ is nonnegative, periodic
in $\R^N$, $p_{\infty}=0$ in $\R^N\backslash\Omega$, the restriction
of $p_{\infty}$ on $\overline{\Omega}$ is of class
$C^{2,\alpha}(\overline{\Omega})$ and solves
\be\label{eqpinfty}\left\{\baa{rcll}
-\nabla\cdot(A(x)\nabla p_{\infty}) & \!\!=\!\! & F(x,p_{\infty}) & \hbox{in }\overline{\Omega},\vspace{3pt}\\
p_{\infty} & \!\!=\!\! & 0 & \hbox{on }\partial\Omega.\eaa\right.
\ee
Lastly, $p_{\infty}\ge p$ in $\overline{\Omega}$, where $p$ is given in Theorem~$\ref{th2}$.
\end{theo}

We point out that, in general, the function $p_{\infty}$ is not
identically equal to the solution~$p$ of~(\ref{statio}) in
$\overline{\Omega}$. However, it is well equal to~$p$
in~$\overline{\Omega}$ if~$F$ fulfills~(\ref{fkppn}) in $\Omega$. We
refer to Remark~\ref{rem31} for more details.

The last result is concerned with the asymptotic behavior as
$n\to+\infty$ of the pulsating travelling fronts of the
type~$\phi_n(x\cdot e-ct,x)$ connecting~$0$ to~$p_n$ (for
problem~(\ref{eq}) in~$\R^N$ with nonlinearities~$f_n$) and of their
minimal speeds~$c^*_n(e)>0$ in any direction~$e$ (when
$\lambda_{1,n}<0$). The limit shall depend strongly on the
direction~$e$ and blocking phenomena may occur in general.

\begin{theo}\label{th4}
Assume that $\lambda_{1,D}<0$ and let $e$ be any given unit vector of $\R^N$.\par
$a)$ The sequence $(c^*_n(e))_{n\in\N}$ is nonincreasing with limit $c^*(e)\ge 0$.
If all connected components $\mathcal{C}$ of $\Omega$ are
bounded in the direction $e$ in the sense that \be\label{discon}
\sup_{x\in\mathcal{C}}|x\cdot e|<+\infty, \ee then $c^*(e)=0$.\par
$b)$ For any $c\ge c^*(e)$ with $c>0$ and for any sequence $(c_n)_{n\in\N}$ such
that $c_n\to c$ as~$n\to+\infty$ and $c_n\ge c^*_n(e)$, the
pulsating travelling fronts $u_n(t,x)=\phi_n(x\cdot e-c_nt,x)$
for~$(\ref{eq})$ in~$\R^N$ with nonlinearity~$f_n$ satisfy
$$u_n(t,x)\to\left\{\baa{ll}u(t,x) & \hbox{in }C^1_t\hbox{ and }C^2_x\hbox{ locally in }\R\times\Omega,\vspace{3pt}\\
0 & \hbox{in }L^1_{loc}(\R\times(\R^N\backslash\Omega))\eaa\right.$$
up to extraction of a subsequence, where $u(t,x)=\phi(x\cdot
e-ct,x)$ is a classical solution of~$(\ref{eqper})$ with $u_t\ge 0$
in~$\R\times\overline{\Omega}$ and~$\phi(s,\cdot)$ is periodic
in~$\overline{\Omega}$ for all $s\in\R$. Moreover, for any
given~$i\in I_-$, one can shift in time the functions $u_n$ so that
$u(-\infty,\cdot)=0$ and~$u(+\infty,\cdot)>0$ in~$\Omega_i$.\par
$c)$ Assume here that $A$ is constant. If there exist a unit vector
$e'\neq\pm e$ and two real numbers $a<b$ such that
\be\label{slab}
\Omega\supset\big\{x\in\R^N,\, a<x\cdot e'<b\big\},
\ee
then $c^*(e)>0$. If there exist a unit vector $e'$, a point $x_0\in\R^N$
and a real number $r>0$ such that~$e'$ is an eigenvector of $A$ with
$e'\cdot e\neq 0$, and \be\label{cylinder}
\Omega\supset\big\{x\in\R^N,\, d(x,x_0+\R e')<r\big\}, \ee where $d$
denotes the Euclidean distance, then $c^*(e)>0$.
\end{theo}

Theorem~\ref{th4} provides a simple geometrical condition for the
blocking of propagation, in a given direction $e$, in the presence
of hostile periodic patches (by blocking, we mean that~$c^*_n(e)\to
0$ as $n\to+\infty$). Consequently, some quantitative estimates of
the spreading speeds of the solutions~$u$ of the Cauchy
problems~(\ref{cauchyn}) with nonlinearities~$f_n$ can be derived.
Indeed, for any compactly supported function~$u_0\not\equiv 0$, the
solution~$u$ of~(\ref{cauchyn}) with nonlinearity~$f_n$ spreads in
the direction $e$ with the spreading speed
$$w^*_n(e)=\min_{\xi\in\mathbb{S}^{N-1},\,\xi\cdot e>0}\frac{c^*_n(e)}{\xi\cdot e},$$
in the sense that $\liminf_{t\to+\infty}u(t,c\,t\,e+x)\ge p_n(x)$
locally uniformly in $x$ if $0\le c<w^*_n(e)$, whereas
$\lim_{t\to+\infty}u(t,c\,t\,e+x)=0$ locally uniformly in $x$ if
$c>w^*_n(e)$ (see \cite{bhna,fg,w}). In particular, $0<w^*_n(e)\le
c^*_n(e)$. Hence, under condition~(\ref{discon}), $c^*(e)=c^*(-e)=0$
and the solution~$u$ of~(\ref{cauchyn}) with nonlinearity~$f_n$
spreads as slowly as wanted in the directions~$\pm e$ when~$n$ is
large enough. In this  case, since all connected components
of~$\Omega$ are bounded in the direction~$e$, pulsating fronts in
the directions~$\pm e$ for problem~(\ref{eqper}) in~$\Omega$ make no
sense even if, under the notation of part~b), the solutions~$u_n$
can be shifted to converge to a non-trivial solution~$u$
of~(\ref{eqper}) in $\R\times\overline{\Omega}$: what happens is
that, in each connected component of~$\Omega_i$, $u$ is just a time
connection between~$0$ and a non-trivial steady state.\par On the
other hand, Theorem~\ref{th4} also gives some simple geometrical
conditions, of the types~(\ref{slab}) or~(\ref{cylinder}), for
non-blocking in the directions~$\pm e$. These conditions mean
that~$\Omega$ contains a slab which is not orthogonal to~$e$, or
contains a cylinder in a direction which is not orthogonal to~$e$.
We do not know however if these conditions are optimal, even
when~$A$ is constant. Lastly, Theorem~\ref{th4} shows the existence
of pulsating fronts for problem~(\ref{eqper}) in~$\Omega$. Assume
for instance that~$\Omega$ is connected, that is $m=1$ under
notation~(\ref{Omegai}). Then, there are pulsating traveling fronts,
in the usual sense, in the direction~$e$, connecting~$0$ to a
non-trivial periodic stationary solution of~(\ref{eqper}).
Furthermore, if~$F$ is of the KPP type~(\ref{fkppn}) in $\Omega$,
the limiting state is unique and is then equal to the function
$p=p_{\infty}$ given in Theorems~\ref{th2} and~\ref{th3} (see
Remark~\ref{rem31} below and the end of the proof of
Theorem~\ref{th4}). However, Theorem~\ref{th4} holds for general
monostable functions~$F$ which may not be of the KPP type and it
gives the first result about the existence of pulsating fronts with
zero Dirichlet boundary condition in periodic domains (which may not be
cylinders).\hfill\break

\noindent{\bf{Outline of the paper.}} Section~\ref{sec2} is devoted
to the proof of Proposition~\ref{pro1} and Theorem~\ref{th2} about
the existence of minimal non-trivial stationary solutions~$p$ of
problems~(\ref{station}) and~(\ref{statio}) respectively, and about
the large-time behavior of the solutions~$u$ of the Cauchy
problems~(\ref{cauchyn}) and~(\ref{cauchy}). Section~\ref{sec3} is
concerned with the proof of Theorem~\ref{th3} and the relationship
between the minimal solutions~$p_n$ of problems~(\ref{station}) with
nonlinearities~$f_n$ and the minimal solution~$p$ of
problem~(\ref{statio}). Lastly, in Section~\ref{sec4}, we do the
proof of Theorem~\ref{th4} and make clear the role of the
geometrical condition~(\ref{discon}) in the blocking process as
$n\to+\infty$.


\SE{Minimal stationary solutions and large-time beha\-vior for the
Cauchy problems~(\ref{cauchyn}) and~(\ref{cauchy})}\label{sec2}

In the first part of this section, we first deal with the elliptic
and parabolic problems~(\ref{station}) and~(\ref{cauchyn}) set in
the whole space~$\R^N$ with the assumption~(\ref{f}) on the
nonlinearity~$f$. Namely, we do the proof of Proposition~\ref{pro1}.
It is based on the elliptic and parabolic maximum principles and on
the construction of suitable subsolutions. Since some parts of the
proof are quite similar to some arguments used in~\cite{bhr1}
and~\cite{bhr3}, they will only be sketched. In the second part of
this section, we will be concerned with the stationary and Cauchy
problems~(\ref{statio}) and~(\ref{cauchy}) posed in the set~$\Omega$
with zero Dirichlet boundary condition. That is, we will do the proof of
Theorem~\ref{th2}, which will itself be inspired by that of
Proposition~\ref{pro1}, but additional difficulties
arise.\hfill\break

\noindent{\bf{Proof of Proposition~\ref{pro1}.}} Let $\varphi$ be
the unique periodic solution of~(\ref{eqvarphi}) such
that~$\max_{\R^N}\!\varphi=1$. Since the principal periodic
eigenvalue $\lambda_1$ of~(\ref{eqvarphi}) is assumed to be
nega\-tive, one can fix $\epsilon_0\in(0,M]$ so that
$f(x,s)\ge\zeta(x)\,s+(\lambda_1/2)s$ for all
$(x,s)\in\R^N\times[0,\epsilon_0]$. Now, for any
$\epsilon\in(0,\epsilon_0]$, there holds \be\label{strictsub}
-\nabla\cdot(A(x)\nabla(\epsilon\varphi))-f(x,\epsilon\varphi)\le-\epsilon\nabla\cdot(A(x)\nabla\varphi)
-\zeta(x)\epsilon\varphi-\frac{\lambda_1}{2}\epsilon\varphi=\frac{\lambda_1}{2}\epsilon\varphi<0
\ee for all $x\in\R^N$. In other words, the functions
$\epsilon\,\varphi$ are strict subsolutions of~(\ref{station}) for
all $\epsilon\in(0,\epsilon_0]$.\par Let now $U$ be the solution of
the Cauchy problem~(\ref{cauchyn}) with initial datum
$U_0=\epsilon_0\varphi$. Since $0<U_0\le M$ and $f(\cdot,M)\le 0$ in
$\R^N$ and since $U_0$ is a subsolution of~(\ref{station}), it
follows that
$$\epsilon_0\varphi(x)\le U(t,x)\le M\hbox{ for all }(t,x)\in\R_+\times\R^N$$
and that $U$ is nondecreasing with respect to $t$. Furthermore, by
uniqueness for the Cauchy problem~(\ref{cauchyn}), $U(t,\cdot)$ is
periodic in~$\R^N$ for each $t\ge 0$. From standard parabolic
estimates, it follows then that
$$U(t,x)\to p(x)\hbox{ as }t\to+\infty\hbox{ uniformly with respect to }x\in\R^N,$$
where $p$ is a $C^{2,\alpha}(\R^N)$ periodic solution
of~(\ref{station}) such that $0<\epsilon_0\varphi=U_0\le p\le
M$.\par Let us then show that $p$ is the minimal positive solution
of~(\ref{station}) (in the class of all positive solutions
of~(\ref{station}), which are not a priori assumed to be periodic).
Let $q$ be any positive solution of~(\ref{station}).
Let~$\lambda_{1,B(y,R),D}$ denote the principal eigenvalue of the
operator
$$-\nabla\cdot(A(x)\nabla)-\zeta(x)$$
in the open Euclidean ball~$B(y,R)$ of center $y\in\R^N$ and
radius~$R>0$, with zero Dirichlet boundary condition on~$\partial
B(y,R)$. For each point~$y\in\R^N$ and~$R>0$, the principal
eigenvalue~$\lambda_{1,B(y,R),D}$ is characterized by the existence
of a function~$\varphi_{y,R}$ of
class~$C^{2,\alpha}(\overline{B(y,R)})$, solving
$$\left\{\baa{rcll}
-\nabla\cdot(A(x)\nabla\varphi_{y,R})-\zeta(x)\varphi_{y,R} & = & \lambda_{1,B(y,R),D}\,\varphi_{y,R}
& \hbox{in }\overline{B(y,R)},\vspace{3pt}\\
\varphi_{y,R} & > & 0 & \hbox{in }B(y,R),\vspace{3pt}\\
\varphi_{y,R} & = & 0 & \hbox{on }\partial B(y,R).\eaa\right.$$ Up
to normalization, one can assume that
$\max_{\overline{B(y,R)}}\varphi_{y,R}=1$, and the functions
$\varphi_{y,R}$ are then unique. As done in~\cite{bhr3}, there holds
$$\lambda_{1,B(y,R),D}\to\lambda_1\hbox{ as }R\to+\infty,$$
uniformly with respect to $y\in\R^N$. Since $\lambda_1<0$, one can
then fix $R>0$ large enough so
that~$\lambda_{1,B(y,R),D}<\lambda_1/2$ for all $y\in\R^N$. Thus,
for each $y\in\R^N$ and $\epsilon\in(0,\epsilon_0]$, the
function~$\epsilon\varphi_{y,R}$ satisfies
\be\label{varphiyR}\baa{rcl}
-\nabla\cdot(A(x)\nabla(\epsilon\varphi_{y,R}))-f(x,\epsilon\varphi_{y,R}) & \!\!\le\!\! &
-\epsilon\nabla\cdot(A(x)\nabla\varphi_{y,R})-\zeta(x)\epsilon\varphi_{y,R}
-\displaystyle\frac{\lambda_1}{2}\epsilon\varphi_{y,R}\vspace{3pt}\\
& \!\!=\!\! & \Big(\lambda_{1,B(y,R),D}-\displaystyle\frac{\lambda_1}{2}\Big)\epsilon\varphi_{y,R}\vspace{3pt}\\
& \!\!<\!\! & 0\eaa \ee in $B(y,R)$. In other words, the functions
$\epsilon\varphi_{y,R}$ are strict subsolutions of~(\ref{station})
in the balls~$B(y,R)$ for all $\epsilon\in(0,\epsilon_0]$. Now, fix
$y\in\R^N$ and observe that~$\min_{\overline{B(y,R)}}q>0$ by
continuity of~$q$. It follows then that
$$\epsilon^*_y:=\sup\Big\{\epsilon\in(0,\epsilon_0],\ \epsilon\varphi_{y,R}\le q\hbox{ in }\overline{B(y,R)}\Big\}$$
is positive. We shall prove that $\epsilon^*_y=\epsilon_0$. Assume
not. Then~$0<\epsilon^*_y<\epsilon_0$
and~$\epsilon^*_y\varphi_{y,R}\le q$ in~$\overline{B(y,R)}$ with
equality somewhere in $\overline{B(y,R)}$. Since $q>0$ and
$\varphi_{y,R}=0$ on $\partial B(y,R)$, the functions
$\epsilon^*_y\varphi_{y,R}$ and~$q$ are equal somewhere at an
interior point, in $B(y,R)$. But $\epsilon^*_y\varphi_{y,R}$ is a
subsolution of~(\ref{station}), from~(\ref{varphiyR}). Since $f$ is
(at least) Lipschitz-continuous locally with respect to the second
variable, uniformly in $x$, it follows from the strong elliptic
maximum principle that~$\epsilon^*_y\varphi_{y,R}=q$ in $B(y,R)$,
which is impossible since the inequality~(\ref{varphiyR}) is strict.
Therefore, $\epsilon^*_y=\epsilon_0$ for all~$y\in\R^N$ and, in
particular,
$$q(y)\ge\epsilon_0\varphi_{y,R}(y)\hbox{ for all }y\in\R^N.$$
But, by uniqueness of the principal eigenfunctions~$\varphi_{y,R}$
and by periodicity of $A$ and $\zeta$, the function
$y\mapsto\varphi_{y,R}(y)$ is continuous and periodic in $\R^N$.
Since it is positive, one gets
that~$\min_{y\in\R^N}\varphi_{y,R}(y)>0$. Therefore,
$\inf_{\R^N}q>0$.\par Define now
$$\epsilon^*=\sup\Big\{\epsilon\in(0,\epsilon_0],\ \epsilon\varphi\le q\hbox{ in }\R^N\Big\},$$
where we recall that $\varphi$ is the unique periodic solution
of~(\ref{eqvarphi}) such that~$\max_{\R^N}\!\varphi=1$. Since~$q$ is
bounded from below in the whole space~$\R^N$ by a positive constant
and since $\varphi$ is bounded, one has $\epsilon^*>0$. Assume that
$\epsilon^*<\epsilon_0$. Then $\epsilon^*\varphi\le q$ in~$\R^N$ and
there exists a sequence~$(x_k)_{k\in\N}$ in~$\R^N$ such that
$$\epsilon^*\varphi(x_k)-q(x_k)\to 0\hbox{ as }k\to+\infty.$$
By writing $x_k=x'_k+x''_k$ with $x'_k\in L_1\Z\times\cdots L_N\Z$
and $x''_k\in[0,L_1]\times\cdots\times[0,L_N]$, it follows that the
functions $q_k(x)=q(x+x'_k)$ converge, up to extraction of a
subsequence, to a solution~$q_{\infty}$ of~(\ref{station}) such that
$\epsilon^*\varphi\le q_{\infty}$ in $\R^N$ with equality somewhere
in~$\R^N$. As above, one concludes that
$\epsilon^*\varphi=q_{\infty}$ in $\R^N$, which is impossible since
$\epsilon^*\varphi$ is a strict subsolution of~(\ref{station}),
from~(\ref{strictsub}). Therefore, $\epsilon^*=\epsilon_0$, whence
$\epsilon_0\varphi\le q$ in $\R^N$. The parabolic maximum principle
implies that
$$U(t,x)\le q(x)\hbox{ for all }(t,x)\in\R_+\times\R^N,$$
where we recall that $U$ denotes the solution of~(\ref{cauchyn})
with initial datum~$\epsilon_0\varphi$. By passing to the limit
as $t\to+\infty$, one gets that
$$p(x)\le q(x)\hbox{ for all }x\in\R^N.$$\par
Finally, let $u_0:\R^N\to[0,M]$ be a uniformly continuous function
which is not identically equal to~$0$, and let $u$ denote the
solution of~(\ref{cauchyn}) with initial datum~$u_0$. The
maximum principle implies that~$0\le u(t,x)\le M$ for all
$(t,x)\in\R_+\times\R^N$, and $u(t,x)>0$ for all $t>0$ and
$x\in\R^N$. With the same notation as above, there exists then
$\epsilon\in(0,\epsilon_0]$ such that
$$\epsilon\,\varphi_{0,R}\le u(1,\cdot)\hbox{ in }\overline{B(0,R)},$$
where we recall that $R>0$ was chosen so
that~$\lambda_{1,B(y,R),D}<\lambda_1/2$ for all $y\in\R^N$. Let $v$
be the solution of~(\ref{cauchyn}) with initial datum
$$v_0(x)=\left\{\baa{ll}
\epsilon\,\varphi_{0,R}(x) & \hbox{if }x\in\overline{B(0,R)},\vspace{3pt}\\
0 & \hbox{if }x\in\R^N\backslash\overline{B(0,R)}.\eaa\right.$$
Since $0\le v_0\le u(1,\cdot)\le M$ in $\R^N$, there holds
$$0\le v(t,x)\le u(t+1,x)\le M\hbox{ for all }(t,x)\in\R_+\times\R^N.$$
Furthermore, since $v_0$ is a subsolution of~(\ref{station}) because
of~(\ref{varphiyR}) and~$v_0=0$ in
$\R^N\backslash\overline{B(0,R)}$, it follows from the maximum
principle that $v$ is nondecreasing with respect to $t$. Hence, from
standard parabolic estimates, one gets that
$$v(t,x)\to v_{\infty}(x)\hbox{ as }t\to+\infty\hbox{ locally uniformly in }x\in\R^N,$$
where $v_{\infty}$ is a solution of~(\ref{station})
satisfying~$v_0\le v_{\infty}\le M$ in $\R^N$. Notice in particular
that~$v_{\infty}$ is positive in $\R^N$ from the strong maximum
principle, since $v_0$ is nonnegative and not identically equal
to~$0$. But the previous paragraphs yield then $v_{\infty}\ge p$.
Therefore,
$$\liminf_{t\to+\infty}u(t,x)\ge p(x)$$
locally uniformly in $x\in\R^N$. Lastly, if $u_0\le p$, then
$u(t,x)\le p(x)$ for all $(t,x)\in\R_+\times\R^N$, whence $u(t,x)\to
p(x)$ as $t\to+\infty$ locally uniformly in $x\in\R^N$. The proof of
Proposition~\ref{pro1} is thereby complete.\hfill$\Box$\break

Let us now turn to the proof of Theorem~\ref{th2}. Some of the ideas
of the proof of Proposition~\ref{pro1} can be adapted. However, the
case of problems~(\ref{statio}) and~(\ref{cauchy}) in $\Omega$ is
substantially more involved than the case of the whole space~$\R^N$,
mainly due to the fact that zero Dirichlet boundary condition is
imposed on~$\partial\Omega$ and the connected components of~$\Omega$
may be bounded or unbounded.\hfill\break

\noindent{\bf{Proof of Theorem~\ref{th2}.}} Remember that the sets
$\Omega_i$ given in~(\ref{Omegai}) are all periodic and pairwise
disjoint. We first work in each set $\Omega_i$ for which
$\lambda_{1,\Omega_i,D}<0$, that is $i\in I_-$. We claim that, for
each such index $i\in I_-$, there exists a periodic solution
$\widetilde{p}_i\in C^{2,\alpha}(\overline{\Omega_i})$ of the
stationary problem \be\label{statioi}\left\{\baa{rcll}
-\nabla\cdot(A(x)\nabla \widetilde{p}_i) & = & F(x,\widetilde{p}_i(x)) & \hbox{in }\overline{\Omega_i},\vspace{3pt}\\
\widetilde{p}_i & = & 0 & \hbox{on }\partial\Omega_i,\vspace{3pt}\\
\widetilde{p}_i & > & 0 & \hbox{in }\Omega_i.\eaa\right. \ee Indeed,
let $\widetilde{\varphi}_i$ be the principal periodic eigenfunction
of the operator $-\nabla\cdot(A(x)\nabla)-\zeta(x)$ in $\Omega_i$
with zero Dirichlet boundary condition on $\partial\Omega_i$. That is,
the function $\widetilde{\varphi}_i$ is periodic, of
class~$C^{2,\alpha}(\overline{\Omega_i})$, and it solves
\be\label{varphii}\left\{\baa{rcll}
-\nabla\cdot(A(x)\nabla\widetilde{\varphi}_i)-\zeta(x)\widetilde{\varphi}_i & = &
\lambda_{1,\Omega_i,D}\,\widetilde{\varphi}_i & \hbox{in }\overline{\Omega_i},\vspace{3pt}\\
\widetilde{\varphi}_i & = & 0 & \hbox{on }\partial\Omega_i,\vspace{3pt}\\
\widetilde{\varphi}_i & > & 0 & \hbox{in }\Omega_i.\eaa\right. \ee
Up to normalization, one can assume that
$\max_{\overline{\Omega_i}}\widetilde{\varphi}_i=1$. Now, as in the
proof of Proposition~\ref{pro1}, since $\lambda_{1,\Omega_i,D}<0$,
there exists $\epsilon_0\in(0,M]$ such that, for any
$\epsilon\in(0,\epsilon_0]$, the
function~$\epsilon\widetilde{\varphi}_i$ is a strict subsolution
of~(\ref{statioi}), namely \be\label{subomegai}
-\nabla\cdot(A(x)\nabla(\epsilon\widetilde{\varphi}_i))-F(x,\epsilon\widetilde{\varphi}_i(x))<0\hbox{
in }\Omega_i, \ee together with $\epsilon\widetilde{\varphi}_i=0$ on
$\partial\Omega_i$ and $\epsilon\widetilde{\varphi}_i>0$ in
$\Omega_i$. But since the constant~$M$ is a supersolution of this
problem, the solution $u_i$ of the Cauchy problem
\be\label{cauchyi}\left\{\baa{rcll}
(u_i)_t-\nabla\cdot(A(x)\nabla u_i) & = & F(x,u_i), & t>0,\ x\in\overline{\Omega_i},\vspace{3pt}\\
u_i(t,x) & = & 0, & t>0,\ x\in\partial\Omega_i,\vspace{3pt}\\
u_i(0,x) & = & \epsilon_0\widetilde{\varphi}_i(x), &
x\in\Omega_i,\eaa\right. \ee is such that
$\epsilon_0\widetilde{\varphi}_i(x)\le u_i(t,x)\le M$ for all
$(t,x)\in(0,+\infty)\times\overline{\Omega_i}$ and $u_i$ is
nondecreasing in~$t$ and periodic in $x$ in $\overline{\Omega_i}$.
Therefore, there exists a
periodic~$C^{2,\alpha}(\overline{\Omega_i})$
solution~$\widetilde{p}_i$ of~(\ref{statioi}) such that $u_i(t,x)\to
\widetilde{p}_i(x)$ as $t\to+\infty$, uniformly in
$x\in\overline{\Omega_i}$.\par Let now $\widetilde{q}_i$ be any
classical bounded solution of~(\ref{statioi}) and let us prove that
$\widetilde{q}_i\ge \widetilde{p}_i$ in $\overline{\Omega_i}$. By
definition of~$\Omega_i$, the set $\omega_i$ is one of its connected
components, and any of its connected components is of the type
$\omega_i+k$ for some $k\in L_1\Z\times\cdots\times L_N\Z$. Two
cases may then occur: either~$\omega_i$ is bounded, or~$\omega_i$ is
unbounded.\par {\it Case 1.} Consider first the case when $\omega_i$
is bounded. Since $\widetilde{q}_i>0$
in~$\omega_i\,(\subset\Omega_i)$, $\widetilde{q}_i=0$
on~$\partial\omega_i\,(\subset\partial\Omega_i)$ and
$F(\cdot,0)\equiv 0$, it follows from Hopf lemma and the compactness
of~$\partial\omega_i$ that
$$\max_{x\in\partial\omega_i}\frac{\partial \widetilde{q}_i}{\partial\nu}(x)<0,$$
where $\nu$ denotes the outward unit normal on $\partial\Omega$. On
the other hand, the principal eigenfunction~$\widetilde{\varphi}_i$
of~(\ref{varphii}) is (at least) of class~$C^1(\overline{\omega_i})$
and $\widetilde{\varphi}_i=0$ on~$\partial\omega_i$. Hence, the
quantity
$$\epsilon^*:=\sup\big\{\epsilon\in(0,\epsilon_0],\ \epsilon\,\widetilde{\varphi}_i\le
\widetilde{q}_i\hbox{ in }\overline{\omega_i}\big\}$$ is a positive
real number, belonging to the interval $(0,\epsilon_0]$.
Furthermore, $\epsilon^*\widetilde{\varphi}_i\le \widetilde{q}_i$
in~$\overline{\omega_i}$. Since $\epsilon^*\widetilde{\varphi}_i$ is
a strict subsolution in $\omega_i\subset\Omega_i$, in the sense
of~(\ref{subomegai}), the strong maximum principle and the Hopf
lemma imply that $\epsilon^*\widetilde{\varphi}_i<\widetilde{q}_i$
in~$\omega_i$ and
$$\frac{\partial \widetilde{q}_i}{\partial\nu}<\epsilon^*\frac{\partial\widetilde{\varphi}_i}{\partial\nu}\
\hbox{ on }\partial\omega_i.$$ Therefore, there exists $\eta_0>0$
such that $(\epsilon^*+\eta)\,\widetilde{\varphi}_i\le
\widetilde{q}_i$ in~$\overline{\omega_i}$ for all
$\eta\in[0,\eta_0]$. The definition of $\epsilon^*$ then yields
$\epsilon^*=\epsilon_0$, whence $\epsilon_0\widetilde{\varphi}_i\le
\widetilde{q}_i$ in $\overline{\omega_i}$. The same argument can be
repeated in~$\omega_i+k$ for all $k\in L_1\Z\times\cdots\times
L_N\Z$. Therefore, $\epsilon_0\widetilde{\varphi}_i\le
\widetilde{q}_i$ in $\overline{\Omega_i}$. By
comparing~$\widetilde{q}_i$ with the solution~$u_i$ of the Cauchy
problem~(\ref{cauchyi}), it follows then as in the proof of
Proposition~\ref{pro1} that \be\label{piqi} \widetilde{p}_i\le
\widetilde{q}_i\ \hbox{ in }\overline{\Omega_i}. \ee\par {\it Case
2.} Consider now the case when $\omega_i$ is unbounded. For all
$y\in\omega_i$ and $R>0$, define
$$\omega_{i,y,R}=\big\{z\in\omega_i,\ d_{\Omega}(y,z)<R\big\},$$
where $d_{\Omega}$ denotes the geodesic distance inside $\Omega$,
and set \be\label{deflambdaomegai}
\lambda_{1,\omega_{i,y,R},D}=\min_{\phi\in
H^1_0(\omega_{i,y,R})\backslash\{0\}}
\frac{\displaystyle{\int_{\omega_{i,y,R}}}A\nabla\phi\cdot\nabla\phi-\zeta\phi^2}
{\displaystyle{\int_{\omega_{i,y,R}}}\phi^2}. \ee Actually,
$\lambda_{1,\omega_{i,y,R},D}$ is the smallest eigenvalue of the
operator $-\nabla\cdot(A\nabla)-\zeta$ in $\omega_{i,y,R}$ with
zero Dirichlet boundary condition (that is, in the
$H^1_0(\omega_{i,y,R})$ sense), but, since $\partial\omega_{i,y,R}$
may not be smooth in general, the
eigenvalue~$\lambda_{1,\omega_{i,y,R},D}$ may not be associated
with~$C^1(\overline{\omega_{i,y,R}})$ eigenfunctions. We first claim
that
$$\limsup_{R\to+\infty}\Big(\sup_{y\in\omega_i}\lambda_{1,\omega_{i,y,R},D}\Big)<0.$$
To do so, let $\widetilde{\rho}:\R\to[0,1]$ be a $C^{\infty}(\R)$
function such that $\widetilde{\rho}=1$ on $(-\infty,-1]$ and
$\widetilde{\rho}=0$ on~$[0,+\infty)$ and, for all $y\in\omega_i$
and $R>0$, denote
$$\rho_{y,R}(x)=\widetilde{\rho}\big(d_{\Omega}(x,y)-R\big)\ \hbox{ for all }x\in\omega_i.$$
These functions $\rho_{y,R}$ are then is $W^{1,\infty}(\omega_i)$.
For every $y\in\omega_i$ and $R>1$, the restriction of the function
$\widetilde{\varphi}_i\rho_{y,R}$ to $\omega_{i,y,R}$ belongs to
$H^1_0(\omega_{i,y,R})\backslash\{0\}$, whence
$$\baa{rcl}
\lambda_{1,\omega_{i,y,R},D} & \le & \frac{\displaystyle{\int_{\omega_{i,y,R}}}A\nabla(\widetilde{\varphi}_i\rho_{y,R})
\cdot\nabla(\widetilde{\varphi}_i\rho_{y,R})-\zeta(\widetilde{\varphi}_i\rho_{y,R})^2}
{\displaystyle{\int_{\omega_{i,y,R}}}(\widetilde{\varphi}_i\rho_{y,R})^2}\vspace{3pt}\\
& \le & \frac{\displaystyle{\int_{\omega_{i,y,R}}}\rho_{y,R}A\nabla
\widetilde{\varphi}_i\cdot\nabla(\widetilde{\varphi}_i\rho_{y,R})-\zeta(\widetilde{\varphi}_i\rho_{y,R})^2}
{\displaystyle{\int_{\omega_{i,y,R}}}(\widetilde{\varphi}_i\rho_{y,R})^2}+
\displaystyle\frac{M\,|\omega_{i,y,R}\backslash\omega_{i,y,R-1}|}
{\displaystyle{\int_{\omega_{i,y,R}}}(\widetilde{\varphi}_i\rho_{y,R})^2},\eaa$$
where
$$M=(1+\|\nabla\widetilde{\varphi}_i\|_{L^{\infty}(\Omega_i)})\times\max_{x\in
\overline{\Omega},\,|\xi|=1,\,|\xi'|=1}(A(x)\xi\cdot\xi')$$ is a
positive constant which does not depend on $y$ or $R$, and
$|\omega_{i,y,R}\backslash\omega_{i,y,R-1}|$ denotes the Lebesgue
measure of $\omega_{i,y,R}\backslash\omega_{i,y,R-1}$. By
integrating by parts, it follows then from~(\ref{varphii}) that
$$\lambda_{1,\omega_{i,y,R},D}\le\lambda_{1,\Omega_i,D}+\displaystyle\frac{2\,M\,|\omega_{i,y,R}\backslash
\omega_{i,y,R-1}|}{\displaystyle{\int_{\omega_{i,y,R}}}(\widetilde{\varphi}_i\rho_{y,R})^2}.$$
Since $\widetilde{\varphi}_i$ is periodic and positive in $\Omega_i$
(and then uniformly away from $0$ in each non-empty set of the type
$$\omega_i^{\delta}:=\big\{x\in\omega_i,\ d(x,\partial\omega_i)>\delta\big\}$$
with $\delta>0$) and since $\Omega$ (and hence $\omega_i$) has a smooth boundary, it follows that
$$\liminf_{R\to+\infty}\Big(\inf_{y\in\omega_i}|\omega_{i,y,R-1}|^{-1}
\int_{\omega_{i,y,R}}(\widetilde{\varphi}_i\rho_{y,R})^2\Big)>0,$$
while
$\limsup_{R\to+\infty}\big(\sup_{y\in\omega_i}|\omega_{i,y,R-1}|^{-1}|\omega_{i,y,R}\backslash
\omega_{i,y,R-1}|\big)=0$. Remember that $\lambda_{1,\Omega_i,D}<0$.
Therefore, there exists $R_0>1$ such that \be\label{lambdaiyR}
\forall\,R\ge R_0,\ \forall\,y\in\omega_i,\ \
\lambda_{1,\omega_{i,y,R},D}<\frac{\lambda_{1,\Omega_i,D}}{2}.
\ee\par Let now $\delta>0$ be any positive constant such that
$\omega_i^{\delta}\neq\emptyset$ and let us show that
$\inf_{\omega_i^{\delta}}\!\widetilde{q}_i>0$. Assume not and let
$\epsilon_i>0$ be such that
$F(x,s)\ge\zeta(x)s+(\lambda_{1,\Omega_i,D}/2)s$ for all
$(x,s)\in\overline{\Omega}\times[0,\epsilon_i]$. There is then a
sequence $(x_n)_{n\in\\N}$ in $\omega_i^{\delta}$ such that
$\widetilde{q}_i(x_n)\to0$ as $n\to+\infty$. Since
$\widetilde{q}_i\ge0$ in $\overline{\omega_i}$ and
$\widetilde{q}_i=0$ on $\partial\omega_i$, it follows from Harnack
inequality that
$$\max_{\overline{\omega_{i,x_{n_0},R_0}}}\widetilde{q}_i\le\epsilon_i\ \hbox{ for some }n_0\in\N\hbox{ large enough}.$$
In particular, \be\label{qiomega} -\nabla\cdot(A(x)\nabla
\widetilde{q}_i)-\zeta(x)\widetilde{q}_i\ge\frac{\lambda_{1,\Omega_i,D}}{2}\,\widetilde{q}_i\
\hbox{ in }\omega_{i,x_{n_0},R_0}. \ee On the other hand, from
(\ref{deflambdaomegai}) and~(\ref{lambdaiyR}), and owing to the
definition of $H^1_0(\omega_{i,x_{n_0},R_0})$, there is $\phi\in
C^1_c(\omega_{i,x_{n_0},R_0})\backslash\{0\}$ (with a compact
support which is included in $\omega_{i,x_{n_0},R_0}$) such that
$$R[\phi]:=\frac{\displaystyle{\int_{\omega_{i,x_{n_0},R_0}}}A\nabla\phi\cdot\nabla\phi-\zeta\phi^2}
{\displaystyle{\int_{\omega_{i,x_{n_0},R_0}}}\phi^2}<\frac{\lambda_{1,\Omega_i,D}}{2}.$$
Now, let $\omega'$ be any bounded open set of class $C^{2,\alpha}$,
containing the support of $\phi$, and such
that~$\overline{\omega'}\subset\omega_{i,x_{n_0},R_0}$. It follows
that $\lambda_{1,\omega',D}\le R[\phi]<\lambda_{1,\Omega_i,D}/2$.
There is then a nonnegative and nontrivial function $\varphi'\in
C^{2,\alpha}(\overline{\omega'})$ solving \be\label{omega'}
-\nabla\cdot(A(x)\nabla\varphi')-\zeta(x)\varphi'=\lambda_{1,\omega',D}\varphi'\le
\frac{\lambda_{1,\Omega_i,D}}{2}\varphi'\ \hbox{ in }\omega' \ee
with $\varphi'=0$ on $\partial\omega'$. Notice that $\varphi'$ may
not be positive in $\omega'$ since $\omega'$ may not be connected.
But $\varphi'$ is positive at least in one connected
component~$\omega''$ of $\omega'$. Since
$\min_{\overline{\omega''}}\widetilde{q}_i>0$ and $\varphi'=0$ on
$\partial\omega''$, it follows from~(\ref{qiomega}),~(\ref{omega'})
and the strong maximum principle that $\epsilon\varphi'\le
\widetilde{q}_i$ in $\overline{\omega''}$ for all $\epsilon>0$,
which is clearly impossible. One has then reached a contradiction.
Hence there holds \be\label{infqi}
\inf_{\omega_i^{\delta}}\widetilde{q}_i>0\ \hbox{ for all
}\delta>0\hbox{ such that }\omega_i^{\delta}\neq\emptyset. \ee\par
It follows then from~(\ref{infqi}), together with the Hopf lemma and
the global smoothness of~$\partial\omega_i$, that
$\sup_{\partial\omega_i}\frac{\partial
\widetilde{q}_i}{\partial\nu}<0$. Therefore, the quantity
$$\epsilon^*:=\sup\big\{\epsilon\in(0,\epsilon_0],\ \epsilon\widetilde{\varphi}_i\le \widetilde{q}_i\hbox{ in }
\overline{\omega_i}\big\}$$ is a positive real number.
From~(\ref{subomegai}) and the strong maximum principle, there holds
$\epsilon^*\widetilde{\varphi}_i<\widetilde{q}_i$ in $\omega_i$.
Furthermore, we claim that \be\label{claim}
\inf_{\omega_i^{\delta}}(\widetilde{q}_i-\epsilon^*\widetilde{\varphi}_i)>0\
\hbox{ for all }\delta>0\hbox{ such that
}\omega_i^{\delta}\neq\emptyset. \ee Assume not. Then there exist
$\delta>0$ such that $\omega_i^{\delta}\neq\emptyset$ and a
sequence~$(y_n)_{n\in\N}$ in $\omega_i^{\delta}$ such
that~$\widetilde{q}_i(y_n)-\epsilon^*\widetilde{\varphi}_i(y_n)\to
0$ as $n\to+\infty$. Write $y_n=y'_n+y''_n$ where $y'_n\in
L_1\Z\times\cdots\times L_N\Z$ and $y''_n\in C$.  Notice in
particular that
$d(y''_n,\partial\Omega)=d(y''_n,\partial\Omega_i)>\delta$. Up to
extraction of a subsequence, one can assume that $y''_n\to
y_{\infty}\in\Omega_i$ as $n\to+\infty$ with
$$d(y_{\infty},\partial\Omega)=d(y_{\infty},\partial\Omega_i)\ge\delta,$$
and that the functions $x\mapsto \widetilde{q}_i(x+y'_n)$ defined in
$\overline{\Omega_i}$ converge in $C^2_{loc}(\overline{\Omega_i})$
to a solution $\overline{q}_i$ of
$$-\nabla\cdot(A(x)\nabla\overline{q}_i)=F(x,\overline{q}_i(x))\hbox{ in }\overline{\Omega_i}$$
such that $\overline{q}_i\ge\epsilon^*\widetilde{\varphi}_i$ in
$B(y_{\infty},\delta)\subset\Omega_i$ with equality at $y_{\infty}$.
The strong maximum principle and~(\ref{subomegai}) lead to a
contradiction. Thus, the claim~(\ref{claim}) holds. As above, it
follows then from Hopf lemma and the global smoothness of
$\partial\omega_i$ that
$\sup_{\partial\omega_i}\frac{\partial(\widetilde{q}_i-\epsilon^*\widetilde{\varphi}_i)}{\partial\nu}<0$
and that there exists $\eta_0>0$ such that
$(\epsilon^*+\eta)\,\widetilde{\varphi}_i\le \widetilde{q}_i$ in
$\overline{\omega_i}$ for all $\eta\in[0,\eta_0]$.
Therefore,~$\epsilon^*=\epsilon_0$,
whence~$\epsilon_0\widetilde{\varphi}_i\le \widetilde{q}_i$ in
$\overline{\omega_i}$ and then in $\overline{\Omega_i}$ by repeating
the argument in $\omega_i+k$ for all $k$ in~$L_1\Z\times\cdots\times
L_N\Z$. Finally, by comparing $\widetilde{q}_i$ with the
solution~$u_i$ of the Cauchy problem~(\ref{cauchyi}), the
conclusion~(\ref{piqi}) follows.\par {\it Conclusion of the proof.}
Define the function $p$ in $\overline{\Omega}$ by
$$p=\left\{\baa{ll}
\widetilde{p}_i & \hbox{in all the sets }\overline{\Omega_i}\hbox{ with }i\in I_-,\vspace{3pt}\\
0 & \hbox{in all the sets }\overline{\Omega_i}\hbox{ with }i\not\in
I_-.\eaa\right.$$ The function $p$ is periodic, of class
$C^{2,\alpha}(\overline{\Omega})$, and it solves~(\ref{statio}).
Furthermore, it follows from the previous steps that any bounded
solution $q$ of~(\ref{statio}) is such that $q\ge p$ in
$\overline{\Omega}$. Lastly, let~$u_0:\overline{\Omega}\to[0,M]$ be
any uniformly continuous function such that $u_0\not\equiv 0$ in
$\overline{\Omega}$, let $u$ be the solution of the Cauchy
problem~(\ref{cauchy}) and let $\omega$ be a connected component of
$\Omega$ intersecting the support of $u_0$. We shall prove that
\be\label{minK} \liminf_{t\to+\infty}\Big(\min_{x\in
K}(u(t,x)-p(x))\Big)\ge 0 \ee for any compact set
$K\subset\overline{\omega}$. Since $0\le u(t,\cdot)\,(\le M)$ in
$\overline{\Omega}$ for all $t>0$ and $p=0$ in $\overline{\Omega_i}$
for all $i\not\in I_i$, it is sufficient to consider the case when
$\omega=\omega_i+k$ for some $i\in I_-$ and some~$k\in
L_1\Z\times\cdots\times L_N\Z$.\par If $\omega$ is bounded, then
$u(1,\cdot)>0$ in $\omega$ and $\max_{\partial\omega}\frac{\partial
u(1,\cdot)}{\partial\nu}<0$ from the strong parabolic maximum
principle. Therefore, $u(1,\cdot)\ge\epsilon\widetilde{\varphi}_i$
in $\overline{\omega}$ for some $\epsilon\in(0,\epsilon_0]$ and
$$u(t+1,x)\ge v(t,x)\hbox{ for all }(t,x)\in(0,+\infty)\times\overline{\omega},$$
where $v$ is the solution of the Cauchy problem~(\ref{cauchyi}) in
$\overline{\omega}$ with initial datum
$\epsilon\widetilde{\varphi}_i$ in~$\overline{\omega}$ and zero Dirichlet
boundary condition on $\partial\omega$. Owing to~(\ref{subomegai}),
$v(t,x)$ is increasing with respect to $t$ (and bounded from above
by the constant~$M$), and it converges as $t\to+\infty$ uniformly in
$\overline{\omega}$ to a solution $w$ of~(\ref{statioi}) in
$\overline{\omega}$ such that $w\ge\epsilon\widetilde{\varphi}_i$ in
$\overline{\omega}$ (whence $w>0$ in $\omega$) and $w=0$ 
on~$\partial\omega$. It follows as in the study of case~1 above that
$w\ge p$ in $\overline{\omega}$, which yields~(\ref{minK}).\par
Consider now the case when $\omega$ is unbounded. Without loss of
generality, up to a translation of the origin, one can assume that
$k=0$ and $\omega=\omega_i$. Choose any point $y_0$ in $\omega$ and,
from~(\ref{lambdaiyR}), let $R_0>0$ be such that
$\lambda_{1,\omega_{i,y_0,R_0},D}<0$. As above, there is then a
function $\phi\in C^1_c(\omega_{i,y_0,R_0})\backslash\{0\}$ such
that
$$R'[\phi]:=\frac{\displaystyle{\int_{\omega_{i,y_0,R_0}}}A\nabla\phi\cdot\nabla\phi-\zeta\phi^2}
{\displaystyle{\int_{\omega_{i,y_0,R_0}}}\phi^2}<0$$ and, if
$\omega'$ is any bounded open set of class $C^{2,\alpha}$ containing
the support of $\phi$ and such
that~$\overline{\omega'}\subset\omega_{i,y_0,R_0}\subset\omega$,
there holds $\lambda_{1,\omega',D}\le R'[\phi]<0$. There is then a
nonnegative and nontrivial function $\varphi'\in
C^{2,\alpha}(\overline{\omega'})$ such that
$$-\nabla\cdot(A(x)\nabla\varphi')-\zeta(x)\varphi'=\lambda_{1,\omega',D}\varphi'\ \hbox{ in }\omega'$$
with $\varphi'=0$ on $\partial\omega'$. Therefore, the function
$\epsilon'\varphi'$  is a subsolution of~(\ref{statioi}) in
$\omega'$ for $\epsilon'>0$ small enough and one can also assume
without loss of generality that $\epsilon'\varphi'\le u(1,\cdot)$ in
the compact set $\overline{\omega'}\subset\omega$. Thus, there holds
$u(t+1,x)\ge v(t,x)$ for all
$(t,x)\in(0,+\infty)\times\overline{\omega}$, where~$v$ is the
solution of the Cauchy problem~(\ref{cauchyi}) in
$\overline{\omega}$ with initial datum $v_0=\epsilon'\varphi'$
in~$\overline{\omega'}$ and $v_0=0$ 
in~$\overline{\omega}\backslash\overline{\omega'}$, and zero Dirichlet
boundary condition on $\partial\omega$. But $v(t,x)$ is increasing
with respect to $t$ and bounded from above by~$M$. It converges
locally uniformly in $\overline{\omega}$ to a solution $w$
of~(\ref{statioi}) in~$\overline{\omega}$, such that~$w\ge v_0$
in~$\overline{\omega}$ (whence $w>0$ in $\omega$ from the strong
maximum principle). One concludes as in case~2 above that $w\ge p$
in $\overline{\omega}$, which leads to~(\ref{minK}).\par Lastly,
observe that, if $u_0\le p$ in $\overline{\Omega}$, then
$u(t,\cdot)\le p$ in $\overline{\Omega}$ for all $t>0$.
Hence,~(\ref{minK}) implies that $u(t,x)\to p(x)$ as $t\to+\infty$
uniformly in any compact subset $K\subset\overline{\omega}$, where
$\omega$ is any connected component of $\Omega$ intersecting the
support of~$u_0$. The proof of Theorem~\ref{th2} is thereby
complete.\hfill$\Box$


\SE{Relationship between the problems~(\ref{station}) with
nonlinearities~$f_n$ and~the problem (\ref{eqpinfty})}\label{sec3}

This section is devoted to the proof of Theorem~\ref{th3}. By using
variational arguments,~$H^1$ a priori estimates and Rellich's
theorem, we prove the monotonicity and the convergence of the
principal periodic eigenvalues of the linearized operators in~$\R^N$
associated with the functions~$f_n$, to that of
problem~(\ref{lambda1D}) with zero Dirichlet boundary condition on
$\partial\Omega$. Then, we show the monotonicity and the convergence
of the functions~$p_n$ to a solution $p_{\infty}\ge p$
of~(\ref{eqpinfty}). The minimality of each solution~$p_n$ and
of~$p$ will also be used.\hfill\break

\noindent{\bf Proof of Theorem~\ref{th3}.} Let, for each $n\in\N$,
$\lambda_{1,n}$ and $\varphi_n$ be the principal eigenvalue and
periodic eigenfunction solving~(\ref{lambda1n}). Let $\lambda_{1,D}$
and $\varphi$ solve~(\ref{lambda1D}), where one can always assume
that $\varphi>0$ in each $\Omega_i$ with
$\lambda_{1,\Omega_i,D}=\lambda_{1,D}$, that is $i\in I_{min}$. Call
$$H^1_{per}(\R^N)=\big\{\phi\in H^1_{loc}(\R^N),\ \phi\hbox{ is periodic}\big\},\quad
L^2_{per}(\R^N)=\big\{\phi\in L^2_{loc}(\R^N),\ \phi\hbox{ is periodic}\big\}$$
and $C_0=[0,L_1]\times\cdots\times[0,L_N]$. For each $n\in\N$, there holds
$$\lambda_{1,n}=\min_{\phi\in H^1_{per}(\R^N)\backslash\{0\}}R_n[\phi]=R_n[\varphi_n],$$
where
$$R_n[\phi]=\frac{\displaystyle{\int_{C_0}}A\nabla\phi\cdot\nabla\phi-\zeta_n\phi^2}
{\displaystyle{\int_{C_0}}\phi^2}.$$ Since the sequence
$(\zeta_n(x))_{n\in\N}$ is nonincreasing for each $x\in\R^N$, it
follows that the sequence~$(\lambda_{1,n})_{n\in\N}$ is
nondecreasing.\par We now claim that $\lambda_{1,n}<\lambda_{1,D}$
for each $n\in\N$. The proof is based on some standard comparison
arguments, used in~\cite{bhr3,bnv}. We just sketch it here for the
sake of completeness. Assume that $\lambda_{1,n}\ge\lambda_{1,D}$
for some $n\in\N$. Pick any index $i\in I_{min}$. Since
$\zeta=\zeta_n$ in $\overline{\Omega_i}\subset\overline{\Omega}$,
there holds
$$-\nabla\cdot(A(x)\nabla\varphi_n)-\zeta(x)\varphi_n=\lambda_{1,n}\varphi_n\ge
\lambda_{1,D}\varphi_n\ \hbox{ in }\overline{\Omega_i}$$ and
$\min_{\overline{\Omega_i}}\varphi_n>0$. In other words, the
periodic function $\varphi_n$ is a supersolution of the linear
equation satisfied by the periodic function~$\varphi$ in~$\Omega_i$.
Since $\varphi=0$ on~$\partial\Omega_i$ and~$\varphi$ is (at least)
of class~$C^2(\overline{\Omega_i})$, it follows from the strong
elliptic maximum principle that the quantity
$$\epsilon^*=\sup\big\{\epsilon\in(0,+\infty),\ \epsilon\varphi\le\varphi_n\hbox{ in }\overline{\Omega_i}\big\}$$
is actually equal to~$+\infty$. This is a contradiction since
$\varphi$ is positive in $\Omega_i$. Therefore,
$\lambda_{1,n}<\lambda_{1,D}$ for all $n\in\N$.\par As a
consequence, the sequence $(\lambda_{1,n})_{n\in\N}$ converges
monotonically to a real number~$\lambda_{1,\infty}$ such that
$\lambda_{1,\infty}\le\lambda_{1,D}$. Let us now show that
$\lambda_{1,\infty}=\lambda_{1,D}$. Normalize here the
eigenfunctions~$\varphi_n$ so that $\|\varphi_n\|_{L^2(C_0)}=1$. It
follows that
$$\int_{C_0}A\nabla\varphi_n\cdot\nabla\varphi_n=\lambda_{1,n}+\int_{C_0}\zeta_n\varphi_n^2\le
\lambda_{1,\infty}+\int_{C_0}\zeta_0\varphi_n^2\le\lambda_{1,\infty}+\max_{\R^N}\zeta_0.$$
Thus, the sequence $(\varphi_n)_{n\in\N}$ is bounded in $H^1(C_0)$.
There exists then a function~$\varphi_{\infty}$ in~$H^1_{per}(\R^N)$
such that, up to extraction of a subsequence,
$\varphi_n\to\varphi_{\infty}$ weakly in $H^1_{per}(\R^N)$ and
strongly in $L^2_{per}(\R^N)$. In particular, $\varphi_{\infty}\ge
0$ a.e. in $\R^N$ and $\|\varphi_{\infty}\|_{L^2(C_0)}=1$. Let $K$
be any compact set such that
$K\subset(\R^N\backslash\overline{\Omega})\cap C_0$. For all
$n\in\N$, one has
$$\baa{rcl}
-\Big(\displaystyle{\mathop{\max}_K}\,\zeta_n\Big)\displaystyle{\int_K}\varphi_n^2\le
-\displaystyle{\int_K}\zeta_n\varphi_n^2 & \!\!=\!\! & \lambda_{1,n}
-\displaystyle{\int_{C_0}}A\nabla\varphi_n\cdot\nabla\varphi_n
+\displaystyle{\int_{C_0\backslash K}}\zeta_n\varphi_n^2\vspace{3pt}\\
& \!\!\le\!\! & \lambda_{1,\infty}+\displaystyle{\int_{C_0\backslash
K}}\zeta_0\varphi_n^2
\le\lambda_{1,\infty}+\displaystyle{\mathop{\sup}_{C_0\backslash
K}}|\zeta_0|,\eaa$$ whence $\|\varphi_n\|_{L^2(K)}\to 0$ as
$n\to+\infty$ from~(\ref{hypfn}). Thus, $\varphi_{\infty}=0$ a.e. in
$K$, and then a.e. in~$\R^N\backslash\Omega$ and
$\|\varphi_{\infty}\|_{L^2(\Omega\cap C_0)}=1$. Furthermore, since
$\varphi_{\infty}\in H^1_{per}(\R^N)$, one gets that the restriction
of $\varphi_{\infty}$ to $\Omega$ belongs to~$H^1_{0,per}(\Omega)$,
that is the space of periodic~$H^1_{loc}(\overline{\Omega})$
functions whose trace is equal to~$0$ on~$\partial\Omega$. Lastly,
observe that
$$\int_{\Omega\cap C_0}A\nabla\varphi_n\cdot\nabla\varphi_n
\le\int_{C_0}A\nabla\varphi_n\cdot\nabla\varphi_n
=\lambda_{1,n}+\int_{C_0}\zeta_n\varphi_n^2\le\lambda_{1,\infty}+\int_{C_0}\zeta_0\varphi_n^2,$$
while
$$\int_{C_0}\zeta_0\varphi_n^2\to\int_{C_0}\zeta_0\varphi_{\infty}^2=\int_{\Omega\cap C_0}\zeta\varphi_{\infty}^2$$
as $n\to+\infty$. Therefore,
$$\int_{\Omega\cap C_0}A\nabla\varphi_{\infty}\cdot\nabla\varphi_{\infty}
\le\liminf_{n\to+\infty}\int_{\Omega\cap
C_0}A\nabla\varphi_n\cdot\nabla\varphi_n\le\lambda_{1,\infty}
+\int_{\Omega\cap C_0}\zeta\varphi_{\infty}^2,$$ that is
$R_{\infty}[\varphi_{\infty}]\le\lambda_{1,\infty}\le\lambda_{1,D}$,
where the functional $R_{\infty}$ is defined by
$$R_{\infty}[\phi]=\frac{\displaystyle{\int_{\Omega\cap C_0}}A\nabla\phi\cdot\nabla\phi-\zeta\phi^2}
{\displaystyle{\int_{\Omega\cap C_0}}\phi^2}$$ for all $\phi\in
H^1_{0,per}(\Omega)\backslash\{0\}$. But since $\min_{\phi\in
H^1_{0,per}(\Omega)\backslash\{0\}}R_{\infty}[\phi]=\lambda_{1,D}$,
one concludes that~$\lambda_{1,\infty}=\lambda_{1,D}$. In other
words, $\lambda_{1,n}\to\lambda_{1,D}$ as $n\to+\infty$.\par In the
sequel, assume now that $\lambda_{1,D}<0$. Consequently, for each
$n\in\N$, one has $\lambda_{1,n}<\lambda_{1,D}<0$ and, from
Proposition~\ref{pro1}, there exists a minimal periodic
solution~$p_n$ of
$$\left\{\baa{rcll}
-\nabla\cdot(A(x)\nabla p_n) & = & f_n(x,p_n) & \hbox{in }\R^N,\vspace{3pt}\\
0 & < & p_n\,\le\,M & \hbox{in }\R^N.\eaa\right.$$
Fix any two integers $n\le m$. Since
$$-\nabla\cdot(A(x)\nabla p_n)-f_m(x,p_n)=f_n(x,p_n)-f_m(x,p_n)\ge 0\hbox{ in }\R^N,$$
the function $p_n$ is a supersolution for the equation satisfied by
$p_m$. From the proof of Proposition~\ref{pro1}, there
exists~$\epsilon_m>0$ such that all functions $\epsilon\varphi_m$
with $\epsilon\in(0,\epsilon_m]$ are subsolutions of~(\ref{station})
with the nonlinearity~$f_m$. Since $\min_{\R^N}\!p_n>0$, there
exists $\epsilon\in(0,\epsilon_m]$ such that~$\epsilon\varphi_m\le
p_n$ in~$\R^N$. Hence, the maximum principle implies that
$$v(t,x)\le p_n(x)\hbox{ for all }(t,x)\in\R_+\times\R^N,$$
where~$v$ is the solution of the Cauchy problem~(\ref{cauchyn}) with
the nonlinearity~$f_m$ and initial datum~$v_0=\epsilon\varphi_m$.
But $v$ is nondecreasing in~$t$ and converges as $t\to+\infty$ to a
solution~$q$ of~(\ref{station}) with nonlinearity $f_m$, such that
$0<q\le p_n$ in $\R^N$. By minimality of~$p_m$ (from
Proposition~\ref{pro1}), one gets that $p_m\le q$ in $\R^N$, whence
$$p_m\le p_n\hbox{ in }\R^N.$$
In other words, the sequence of functions~$(p_n)_{n\in\N}$ is
nonincreasing and then converges pointwise to a periodic function
$p_{\infty}(x)$ ranging in $[0,M]$.\par Let us now show that
$p_{\infty}=0$ in $\R^N\backslash\overline{\Omega}$. By multiplying
by~$p_n$ the equation~(\ref{station}) with the nonlinearity~$f_n$,
that is $-\nabla\cdot(A(x)\nabla p_n)=f_n(x,p_n)$, and by
integrating by parts over the cell~$C_0$, it follows that
$$\int_{C_0}A\nabla p_n\cdot\nabla p_n=\int_{C_0}f_n(x,p_n)\,p_n\le\int_{C_0}f_0(x,p_n)\,p_n
\le M\times\max_{\R^N\times[0,M]}|f_0|,$$
whence the sequence $(p_n)_{n\in\N}$ is bounded in
$H^1_{per}(\R^N)$. Since it converges monotonically to~$p_{\infty}$,
one infers that $p_{\infty}\in H^1_{per}(\R^N)$ and $p_n\to
p_{\infty}$ as $n\to+\infty$ weakly in $H^1_{per}(\R^N)$ and
strongly in $L^2_{per}(\R^N)$. For any compact set~$K$ such that
$K\subset(\R^N\backslash\overline{\Omega})\cap C_0$, one has
$$\baa{rcl}
-\Big(\displaystyle{\mathop{\max}_{K\times[0,M]}}g_n\Big)\displaystyle{\int_K}p_n^2
\le-\displaystyle{\int_K}g_n(x,p_n)\,p_n^2 & \!\!=\!\! &-\displaystyle{\int_K}f_n(x,p_n)p_n\vspace{3pt}\\
& \!\!=\!\! & -\displaystyle{\int_{C_0}}A\nabla p_n\cdot\nabla p_n
+\displaystyle{\int_{C_0\backslash K}}f_n(x,p_n)\,p_n\vspace{3pt}\\
& \!\!\le\!\! & \displaystyle{\int_{C_0\backslash K}}f_0(x,p_n)\,p_n
\le
M\times\displaystyle{\mathop{\max}_{\R^N\times[0,M]}}|f_0|.\eaa$$
The assumption~(\ref{hypfn}) yields
$\max_{K\times[0,M]}g_n\to-\infty$ as $n\to+\infty$, whence
$p_{\infty}=0$ a.e. in any such compact~$K$. Finally,~$p_{\infty}=0$
a.e. in $\R^N\backslash\Omega$. Therefore, the restriction
of~$p_{\infty}$ on~$\Omega$ is in~$H^1_{0,per}(\Omega)$.
Furthermore, since \be\label{eqpnOmega} -\nabla\cdot(A(x)\nabla
p_n)=F(x,p_n)\hbox{ in }\overline{\Omega}, \ee the function
$p_{\infty}$ is a solution of the same equation in $\Omega$ in the
weak $H^1_{0,per}(\Omega)$ sense. The elliptic regularity theory
then implies that, up a negligible set,~$p_{\infty}$ is actually a
$C^{2,\alpha}(\overline{\Omega})$ solution of~(\ref{eqpinfty}) and
the convergence $p_n\to p_{\infty}$ holds at least in
the~$C^2_{loc}(\Omega)$ sense.\par Lastly, let us show that
$p_{\infty}\ge p$ in $\Omega$. Since $p_{\infty}$ is nonnegative and
$p=0$ in all~$\Omega_i$ with~$i\not\in I_-$, one only needs to prove
that $p_n\ge p$ in $\Omega_i$ for all $i\in I_-$ and for all
$n\in\N$. For any~$n\in\N$ and~$i\in I_-$, observe that the function
$p_n$ is a supersolution of~(\ref{statioi}) in $\Omega_i$, because
it solves~(\ref{eqpnOmega}) in $\overline{\Omega_i}$ and $p_n>0$ on
$\partial\Omega_i$. Since
$\min_{\overline{\Omega_i}}p_n\ge\min_{\R^N}p_n>0$, there is
$\epsilon>0$ such that~$\epsilon\widetilde{\varphi}_i\le p_n$
in~$\overline{\Omega_i}$, where~$\widetilde{\varphi}_i$
solves~(\ref{varphii}). Since $\lambda_{1,\Omega_i,D}<0$, one can
even assume without loss of generality that
$\epsilon\widetilde{\varphi}_i$ is a subsolution of~(\ref{statioi}),
in the sense of~(\ref{subomegai}). Therefore,
$$w(t,x)\le p_n(x)\hbox{ for all }(t,x)\in\R_+\times\overline{\Omega_i},$$
where~$w$ denotes the solution of the Cauchy problem~(\ref{cauchyi})
in~$\overline{\Omega_i}$ with initial
datum~$\epsilon\widetilde{\varphi}_i$. Since $w$ is nondecreasing in
$t$, it converges as $t\to+\infty$ to a solution~$w_{\infty}$
of~(\ref{statioi}) such that $0<\epsilon\widetilde{\varphi}_i\le
w_{\infty}\le p_n$ in $\Omega_i$. From the construction of~$p$ in
Theorem~\ref{th2} and its minimality, one infers that~$p\le
w_{\infty}$ in $\overline{\Omega_i}$, whence
$$p\le p_n\hbox{ in }\overline{\Omega_i}.$$
As a conclusion, $p\le p_n$ in $\overline{\Omega}$ for all $n\in\N$,
whence $p\le p_{\infty}$ in $\overline{\Omega}$. The proof of
Theorem~\ref{th3} is thereby complete.\hfill$\Box$

\begin{rem}\label{rem31}{\rm We first show in this remark that, if $F$ fulfills the KPP condition~$(\ref{fkppn})$
in $\Omega$, that is if \be\label{Fkpp}
s\mapsto\frac{F(x,s)}{s}\hbox{ is decreasing in }s>0\hbox{ for all
}x\in\Omega, \ee then $p_{\infty}=p$ in $\overline{\Omega}$.
Consider first $i\in I_-$ and let us prove that the function $p$
solving~$(\ref{statioi})$ in~$\overline{\Omega_i}$ is unique. The
proof is similar to the ones used for instance in~\cite{b,bhr1,bhr3}
and it is just sketched. Let $q$ be any periodic solution
of~$(\ref{statioi})$ in~$\overline{\Omega_i}$. From the proof of
Theorem~$\ref{th2}$, one knows that $q\ge p$ in
$\overline{\Omega_i}$. But $\epsilon q\le p$
in~$\overline{\Omega_i}$ for $\epsilon>0$ small enough, from the
Hopf lemma applied to $p$. Therefore, the quantity
$\epsilon^*=\sup\big\{\epsilon>0,\ \epsilon q\le p\hbox{ in
}\overline{\Omega_i}\big\}$ is a positive real number. If
$\epsilon^*<1$, then
$$-\nabla\cdot(A(x)\nabla(\epsilon^*q))-F(x,\epsilon^*q)<0\hbox{ in }\Omega_i,$$
from~$(\ref{Fkpp})$. The strong maximum principle and Hopf lemma
then imply that~$\epsilon^*q<p$ in~$\Omega_i$ and even
$(\epsilon^*+\eta)q<p$ in~$\Omega_i$ for all $\eta\in[0,\eta_0]$ and
for some $\eta_0>0$. This contradicts the maximality of
$\epsilon^*$. Consequently, $\epsilon^*\ge 1$, whence $q\le p$ in
$\overline{\Omega_i}$ and finally $q=p$ in $\overline{\Omega_i}$.
Actually, with the same arguments as those used in the proof of
Theorem~$\ref{th2}$, the same conclusion holds even if $q$ is not
assumed to be periodic. Now, if $i\not\in I_-$, then we prove that
there does not exist any solution $q$ of~$(\ref{statioi})$ that is
positive in~$\Omega_i$ (or in any of its connected components).
Indeed, since $F(x,s)<\zeta(x)s$ for all $x\in\Omega_i$ and $s>0$
from~$(\ref{Fkpp})$, there holds
$$-\nabla\cdot(A(x)\nabla(\epsilon\widetilde{\varphi}_i))-F(x,\epsilon\widetilde{\varphi}_i)>
-\epsilon\nabla\cdot(A(x)\nabla\widetilde{\varphi}_i)-\zeta(x)\epsilon\widetilde{\varphi}_i
=\lambda_{1,\Omega_i,D}\epsilon\widetilde{\varphi}_i\ge 0\hbox{ in
}\Omega_i$$ for all $\epsilon>0$, where $\widetilde{\varphi}_i$
solves~$(\ref{varphii})$ in $\Omega_i$, with
$\lambda_{1,\Omega_i,D}\ge 0$. In other words,
$\epsilon\widetilde{\varphi}_i$ is a strict supersolution
of~$(\ref{statioi})$ for all $\epsilon>0$. It follows with the same
arguments as above or as in the proof of Theorem~$\ref{th2}$ that
$q\le\epsilon\widetilde{\varphi}_i$ for all $\epsilon>0$, for any
solution~$q$ of~$(\ref{statioi})$. Therefore, a positive periodic
solution of~$(\ref{statioi})$ cannot exist, which implies that
$p_{\infty}=p=0$ in $\overline{\Omega_i}$ for all~$i\not\in I_-$. As
a conclusion, the condition~$(\ref{Fkpp})$ implies that
\be\label{pinftyp} p_{\infty}=p\hbox{ in }\overline{\Omega}. \ee

On the other hand, we can construct examples for
which~$(\ref{pinftyp})$ does not hold. It is indeed possible to
construct a situation for which $\lambda_{1,D}<0$ and there exist an
index~$j\in\{1,\ldots,m\}$ and $s_0\in(0,M)$ such that $F(x,s)=\lambda s+s^2$ 
for all $x\in\overline{\Omega_j}$ and $s\in[0,s_0]$, where $\lambda>0$ 
denotes the principal periodic eigenvalue of the
operator $-\nabla\cdot(A(x)\nabla)$ in $\Omega_j$ with zero
Dirichlet boundary condition on $\partial\Omega_j$. Thus,
$\lambda_{1,\Omega_j,D}=0$ and $j\not\in I_-$. Let
$\widetilde{\varphi}_j$ be the principal periodic eigenfunction
of~$(\ref{varphii})$ in $\overline{\Omega_j}$ with $\zeta=\lambda$
in $\overline{\Omega_j}$, such that
$\max_{\overline{\Omega_j}}\widetilde{\varphi}_j=1$. For any
$\epsilon\in(0,s_0]$, there holds
$$-\nabla\cdot(A(x)\nabla(\epsilon\widetilde{\varphi}_j))-F(x,\epsilon\widetilde{\varphi}_j)
=-\epsilon\nabla\cdot(A(x)\nabla\widetilde{\varphi}_j)-\lambda\epsilon\widetilde{\varphi}_j
-\epsilon^2\widetilde{\varphi}_j^2=-\epsilon^2\widetilde{\varphi}_j^2<0\hbox{
in }\Omega_j.$$ As above, it follows from the strong maximum
principle that $\epsilon\widetilde{\varphi}_j\le p_n$ in
$\overline{\Omega_j}$ for all $\epsilon\in(0,s_0]$ and for all
$n\in\N$. In particular, $0<s_0\widetilde{\varphi}_j\le p_{\infty}$
in $\Omega_j$, whereas $p=0$ in $\Omega_j$ by definition.}
\end{rem}


\SE{Pulsating travelling fronts and limiting minimal speed}\label{sec4}

In this section, we give the proof of Theorem~\ref{th4}. We
establish the relationship between the pulsating travelling fronts
for the problems~(\ref{eq}) in~$\R^N$ and~(\ref{eqper}) in~$\Omega$
when the nonlinearity~$F$ is approximated with nonlinearities~$f_n$
which are very negative in~$\R^N\backslash\overline{\Omega}$, in the
sense of~(\ref{hypfn}). We also prove that the minimal speeds of the
fronts in~$\R^N$ converge monotonically to a quantity which is equal
to~$0$ in a direction~$e$ when the connected components of~$\Omega$
are bounded with respect to~$e$. We use especially some bounds for
the minimal speeds, which involve some linear eigenvalue
problems.\par Throughout this section, we assume that
$\lambda_{1,D}<0$ and~$e$ is any given unit vector of~$\R^N$. The
functions~$F$ and~$f_n$ are assumed to fulfill~(\ref{ff})
and~(\ref{hypfn}). For each $n\in\N$, one has~$\lambda_{1,n}<0$ from
Theorem~\ref{th3}. The functions~$p_n$ denote the minimal solutions
of~(\ref{station}) with the nonlinearities~$f_n$, given by
Proposition~\ref{pro1}, and the speeds~$c^*_n(e)>0$ denote the
minimal speeds of pulsating fronts~$\phi_n(x\cdot e-ct,x)$
connecting~$0$ to~$p_n$ for problems~(\ref{eq}) in~$\R^N$ with the
nonlinearities~$f_n$.\hfill\break

\noindent{\bf{Proof of part a) of Theorem~\ref{th4}.}} Fix any two
integers $n\le m$ and let us show that~$c^*_m(e)\le c^*_n(e)$.
First, remember that $0<p_m\le p_n\le M$ (from
Proposition~\ref{pro1} and Theorem~\ref{th3}) and that both
functions~$p_m$ and $p_n$ are periodic in~$\R^N$. Let $\eta>0$ be
such that~$0<\eta<\min_{\R^N}p_m$ and $u_0:\R\to[0,M]$ be defined by
$$u_0(x)=\left\{\baa{ll} 0 & \hbox{if }x\cdot e>0,\vspace{3pt}\\ \eta & \hbox{if }x\cdot e\le 0.\eaa\right.$$
Let $v_n$ and $v_m$ denote the solutions of the Cauchy
problems~(\ref{cauchyn}) with initial datum~$u_0$ and
nonlinearities~$f_n$ and~$f_m$ respectively. Since $f_m\le f_n$, the
maximum principle yields
$$0<v_m(t,x)\le v_n(t,x)<M\ \hbox{ for all }t>0\hbox{ and }x\in\R^N.$$
On the other hand, it follows from the results of Weinberger~\cite{w} that
$$\forall\, c<c^*_m(e),\ \sup_{x\in\R^N,\,x\cdot e\le ct}|v_m(t,x)-p_m(x)|\to 0\hbox{ as }t\to+\infty,$$
while
$$\forall\, c>c^*_n(e),\ \sup_{x\in\R^N,\,x\cdot e\ge ct}v_n(t,x)\to 0\hbox{ as }t\to+\infty.$$
One infers that $c^*_m(e)\le c^*_n(e)$. Consequently, the
sequence~$(c^*_n(e))_{n\in\N}$ is nonincreasing and it converges to
a real number $c^*(e)\ge 0$.\par From the assumptions~(\ref{ff})
and~(\ref{hypfn}) and the regularity of~$F$ and~$f_n$, there exist a
function $\overline{F}:(x,u)\mapsto\overline{F}(x,u)$ and a sequence
of functions~$(\overline{f}_n)_{n\in\N}$ such that: $i$) the
function~$\overline{F}$ is defined and continuous in
$\overline{\Omega}\times\R_+$, of class~$C^{0,\alpha}$ with respect
to~$x\in\overline{\Omega}$ locally uniformly in $u\in\R_+$, of
class~$C^1$ with respect to~$u$ with
$\overline{\zeta}:=\frac{\partial\overline{F}}{\partial
u}(\cdot,0)\in C^{0,\alpha}(\overline{\Omega})$, periodic with
respect to $x\in\overline{\Omega}$ and~$\overline{F}$
satisfies~(\ref{ff}); $ii$) each function $\overline{f}_n$ is
defined and continuous in~$\R^N\times\R_+$, of class~$C^{0,\alpha}$
with respect to $x\in\R^N$ locally uniformly in $u\in\R_+$, of class
$C^1$ with respect to~$u$ with
$\overline{\zeta}_n:=\frac{\partial\overline{f}_n}{\partial
u}(\cdot,0)\in C^{0,\alpha}(\R^N)$, periodic with respect to
$x\in\R^N$ and~$\overline{f}_n$ satisfies~(\ref{f});~$iii$)~the
functions~$\overline{f}_n$ satisfy~(\ref{hypfn}) with $\overline{F}$
instead of~$F$ and $\overline{g}_n(x,u)=\overline{f}_n(x,u)/u$ if
$u>0$, $\overline{g}_n(x,0)=\overline{\zeta}_n(x)$; $iv$) the
function~$\overline{F}$ satisfies
$$F(x,u)\le\overline{F}(x,u)\hbox{ for all }(x,u)\in\overline{\Omega}\times\R_+$$
and $\overline{F}(x,u)/u$ is decreasing with respect to $u>0$ for
all $x\in\Omega$; $v$) the functions~$\overline{f}_n$ satisfy
$$f_n(x,u)\le\overline{f}_n(x,u)\hbox{ for all }(x,u)\in\R^N\times\R_+$$
and $\overline{f}_n(x,u)/u$ is decreasing with respect to $u>0$ for
all $x\in\R^N$.\par Let~$\overline{\lambda}_{1,n}$
and~$\overline{\lambda}_{1,D}$ be the principal periodic eigenvalues
of problems~(\ref{lambda1n}) and~(\ref{lambda1D}) with
coefficients~$\overline{\zeta}_n$ and~$\overline{\zeta}$ instead
of~$\zeta_n$ and~$\zeta$, respectively. Since
$\overline{\zeta}_n\ge\zeta_n$ in~$\R^N$ and
$\overline{\zeta}\ge\zeta$ in~$\overline{\Omega}$, there holds
$$\overline{\lambda}_{1,n}\le\lambda_{1,n}\ \hbox{ and }\ \overline{\lambda}_{1,D}\le\lambda_{1,D},$$
while $\overline{\lambda}_{1,n}<\overline{\lambda}_{1,D}$ and
$\overline{\lambda}_{1,n}\to\overline{\lambda}_{1,D}$ as
$n\to+\infty$ monotonically, from Theorem~\ref{th3}. In particular,
$\overline{\lambda}_{1,n}<\overline{\lambda}_{1,D}<0$ for all
$n\in\N$. Let~$\overline{p}_n$ be the minimal periodic solution
of~(\ref{station}) with the nonlinearity~$\overline{f}_n$, given by
Proposition~\ref{pro1}. Actually, the function~$\overline{p}_n$ is
unique from property~$v$) above and from~\cite{bhr1}, and it is such
that $\overline{p}_n\ge p_n$ in $\R^N$ since $\overline{f}_n\ge f_n$
in $\R^N\times\R_+$, from the proof of Theorem~\ref{th3}. Let
$\overline{c}^*_n(e)>0$ be the minimal speed of pulsating travelling
fronts $\overline{\phi}_n(x\cdot e-ct,x)$ connecting~$0$
to~$\overline{p}_n$ for problem~(\ref{eq}) with the
nonlinearity~$\overline{f}_n$, that is~$\overline{\phi}_n$ is
periodic with respect to~$x\in\R^N$,
$0<\overline{\phi}_n(s,x)<\overline{p}_n(x)$ and
$\overline{\phi}_n(-\infty,x)=\overline{p}_n(x)$,
$\overline{\phi}_n(+\infty,x)=0$. As in the beginning of the proof
of this theorem, there holds \be\label{ccn}
0<c^*_n(e)\le\overline{c}^*_n(e), \ee since $f_n\le\overline{f}_n$.
Furthermore, it follows from~\cite{bhr2,w}
that~$\overline{c}^*_n(e)$ is given by \be\label{c*n}
\overline{c}^*_n(e)=\min_{\lambda>0}\frac{-\overline{k}_{e,\lambda,n}}{\lambda},
\ee where $\overline{k}_{e,\lambda,n}$ denotes the principal
periodic eigenvalue of the operator
$$\overline{\mathcal{L}}_{e,\lambda,n}:=-\nabla\cdot(A\nabla)+2\lambda\,Ae\cdot\nabla+\lambda\,
\nabla\cdot(Ae)-\lambda^2\,Ae\cdot e-\overline{\zeta}_n\ \hbox{ in
}\R^N.$$\par Let us now show that, for every $\lambda\in\R$, one has
$\overline{k}_{e,\lambda,n} \to\overline{k}_{e,\lambda,D}$ as
$n\to+\infty$, where $\overline{k}_{e,\lambda,D}$ is the principal
periodic eigenvalue of the operator
$$\overline{\mathcal{L}}_{e,\lambda,\Omega}:=-\nabla\cdot(A\nabla)+2\lambda\,Ae\cdot\nabla+\lambda\,
\nabla\cdot(Ae)-\lambda^2\,Ae\cdot e-\overline{\zeta}\ \hbox{ in
}\Omega$$ with zero Dirichlet boundary condition on~$\partial\Omega$.
The proof starts as in the proof of the
convergence~$\lambda_{1,n}\to\lambda_{1,D}$ in Theorem~\ref{th3}.
First, it follows as in the proof of Theorem~\ref{th3} that the
sequence~$(\overline{k}_{e,\lambda,n})_{n\in\N}$ is nondecreasing
and that $\overline{k}_{e,\lambda,n}<\overline{k}_{e,\lambda,D}$ for
all $n\in\N$. Let~$\overline{\varphi}_n$ be a principal periodic
eigenfunction of~$\overline{\mathcal{L}}_{e,\lambda,n}$, that is
$$\overline{\mathcal{L}}_{e,\lambda,n}\overline{\varphi}_n=\overline{k}_{e,\lambda,n}\overline{\varphi}_n
\hbox{ and }\overline{\varphi}_n>0\hbox{ in }\R^N.$$ Up to
normalization, one can assume that
$\|\overline{\varphi}_n\|_{L^2(C_0)}=1$. By multiplying the above
equation by~$\overline{\varphi}_n$, by integrating by parts
over~$C_0$ and by using Young's inequality, it follows that the
sequence~$(\overline{\varphi}_n)_{n\in\N}$ is bounded
in~$H^1_{per}(\R^N)$. Up to extraction of a subsequence, it
converges weakly in~$H^1_{per}(\R^N)$ and strongly
in~$L^2_{per}(\R^N)$ to a nonnegative
function~$\overline{\varphi}_{\infty}\in H^1_{per}(\R^N)$ such that
$\|\overline{\varphi}_{\infty}\|_{L^2(C_0)}=1$. Furthermore, since
the sequence $(\overline{k}_{e,\lambda,n})_{n\in\N}$ is bounded and
$\overline{\zeta}_n\to-\infty$ as $n\to+\infty$ locally uniformly in
$\R^N\backslash\overline{\Omega}$, one infers as in the proof of
Theorem~\ref{th3} that~$\overline{\varphi}_{\infty}=0$ a.e. in
$\R^N\backslash\Omega$. The restriction
of~$\overline{\varphi}_{\infty}$ to $\overline{\Omega}$ is then
a~$C^{2,\alpha}(\overline{\Omega})$ periodic function such that
$$\overline{\mathcal{L}}_{e,\lambda,\Omega}\overline{\varphi}_{\infty}
=\overline{k}_{e,\lambda,\infty}\overline{\varphi}_{\infty}\hbox{ in
} \overline{\Omega}\hbox{ with }\overline{\varphi}_{\infty}=0\hbox{
on }\partial\Omega,$$ where
$\lim_{n\to+\infty}\overline{k}_{e,\lambda,n}=\overline{k}_{e,\lambda,\infty}\le
\overline{k}_{e,\lambda,D}$. Since the function
$\overline{\varphi}_{\infty}$ is periodic, nonnegative and
nontrivial, it follows that it is positive in~$\Omega_i$ for some
$i\in\{1,\ldots,m\}$, that is~$\overline{k}_{e,\lambda,\infty}$ is
equal to the principal periodic eigenvalue
$\overline{k}_{e,\lambda,\Omega_i,D}$ of the
operator~$\overline{\mathcal{L}}_{e,\lambda,\Omega}$ in~$\Omega_i$
with zero Dirichlet boundary condition on~$\partial\Omega_i$. But since
$\overline{k}_{e,\lambda,\Omega_i,D}\ge\overline{k}_{e,\lambda,D}\,(\ge\overline{k}_{e,\lambda,\infty})$,
one concludes eventually that
$\overline{k}_{e,\lambda,\infty}=\overline{k}_{e,\lambda,D}$, that
is \be\label{kelD}
\overline{k}_{e,\lambda,n}\to\overline{k}_{e,\lambda,D}\hbox{ as
}n\to+\infty. \ee\par Assume now that all connected components of
$\Omega$ are bounded in the direction~$e$, in the sense
of~(\ref{discon}). Let us show that $c^*(e)=0$. First, it follows
from~(\ref{ccn}),~(\ref{c*n}) and~(\ref{kelD}) that \be\label{c*e}
0\le
c^*(e)\le\inf_{\lambda>0}\frac{-\overline{k}_{e,\lambda,D}}{\lambda}.
\ee On the other hand, for every $\lambda>0$, there is an index
$i\in\{1,\ldots,m\}$, which may depend on~$\lambda$, such that
$\overline{k}_{e,\lambda,D}=\overline{k}_{e,\lambda,\Omega_i,D}$ and
thus there is a periodic function~$\varphi$ defined
in~$\overline{\Omega_i}$ such
that~$\overline{\mathcal{L}}_{e,\lambda,\Omega}\varphi=\overline{k}_{e,\lambda,D}\varphi$
in $\overline{\Omega_i}$ with $\varphi>0$ in~$\Omega_i$ and
$\varphi=0$ on~$\partial\Omega_i$. The function
$\psi=e^{-\lambda(x\cdot e)}\varphi$ satisfies \be\label{eqpsi}
-\nabla\cdot(A(x)\nabla\psi)-\overline{\zeta}(x)\psi=\overline{k}_{e,\lambda,D}\psi\
\hbox{ in }\overline{\Omega_i} \ee with $\psi>0$ in $\Omega_i$ and
$\psi=0$ on $\partial\Omega_i$. Let~$\mathcal{C}$ be any connected
component of~$\Omega_i$, that is~$\mathcal{C}=\omega_i+k$ for some
$k\in L_1\Z\times\cdots\times L_N\Z$. The function~$\psi$ is
positive and bounded in~$\mathcal{C}$ because of~(\ref{discon}) and
since~$\varphi$ is bounded. It follows then from Hopf lemma and the
smoothness of~$\partial\mathcal{C}$ that there exist $r>0$ and a
sequence~$(x_n)_{n\in\N}$ in~$\mathcal{C}$ such that
$B(x_n,r)\subset C$ for all $n\in\N$ and
$\psi(x_n)\to\sup_{\mathcal{C}}\psi$ as $n\to+\infty$. By using the
standard elliptic estimates and passing to the limit
in~(\ref{eqpsi}) in~$B(x_n,r)$, up to extraction of a subsequence,
one infers
that~$\overline{k}_{e,\lambda,D}\ge\liminf_{n\to+\infty}-\overline{\zeta}(x_n)\ge-\max_{\overline{\Omega_i}}
\overline{\zeta}$.
Finally,
$\overline{k}_{e,\lambda,D}\ge-\max_{\overline{\Omega}}\overline{\zeta}$
for all $\lambda>0$, whence $c^*(e)=0$ from~(\ref{c*e}).\hfill\break

\noindent{\bf{Proof of part b) of Theorem~\ref{th4}.}} Let $c$ be
any positive real number such that $c\ge c^*(e)$ and
let~$(c_n)_{n\in\N}$ be any sequence such that $c_n\to c$
as~$n\to+\infty$ and $c_n\ge c^*_n(e)$ for all $n\in\N$. Let
$$u_n(t,x)=\phi_n(x\cdot e-c_nt,x)$$
be pulsating travelling fronts for~$(\ref{eq})$ in~$\R^N$ with nonlinearity~$f_n$, such that
$$0=\phi_n(+\infty,x)<\phi_n(s,x)<\phi_n(-\infty,x)=p_n(x)\le M\ \hbox{ for all }(s,x)\in\R\times\R^N.$$
Actually, from~\cite{h8}, each solution~$u_n$ satisfies $(u_n)_t>0$
in~$\R\times\R^N$.\par On the one hand, since $0<u_n(t,x)<p_n(x)$
in~$\R\times\R^N$, Theorem~\ref{th3} implies that $u_n\to 0$ in
$L^1_{loc}(\R\times(\R^N\backslash\Omega))$. On the other hand,
since $f_n(x,s)=F(x,s)$ for all
$(x,s)\in\overline{\Omega}\times\R_+$, it follows from standard
parabolic estimates that there exists a
function~$u:\R\times\Omega\to[0,M]$ such that, up to extraction of a
subsequence, $u_n\to u$ as~$n\to+\infty$ in~$C^1_t$ and~$C^2_x$
in~$\R\times\Omega$, where~$u$ obeys
$$u_t-\nabla\cdot(A(x)\nabla u)=F(x,u)\ \hbox{ in }\R\times\Omega$$
and $0\le u(t,x)\le p_{\infty}(x)\le M$ for all
$(t,x)\in\R\times\Omega$, under the notation of Theorem~\ref{th3}.
In particular, the function $u$ can be extended continuously by~$0$
on~$\R\times\partial\Omega$ and, from parabolic regularity, the
function $u$ is a classical solution of~(\ref{eqper})
in~$\R\times\overline{\Omega}$ (of course, one could also extend~$u$
by~$0$ in~$\R\times(\R^N\backslash\Omega)$ and $u$ would then be
continuous in $\R\times\R^N$). Moreover, the equalities
$$u_n\Big(t+\frac{k\cdot e}{c_n},x\Big)=u_n(t,x-k)\hbox{ in }\R\times\R^N$$
for all $k\in L_1\Z\times\cdots\times L_N\Z$ carry over at the
limit, whence $u(t+(k\cdot e)/c,x) =u(t,x-k)$ in
$\R\times\overline{\Omega}$ for all $k\in L_1\Z\times\cdots\times
L_N\Z$. In other words, the function~$u$ can be written
as~$u(t,x)=\phi(x\cdot e-ct,x)$ in $\R\times\overline{\Omega}$ where
$\phi:\R\times\overline{\Omega}\to[0,M]$ is such that
$\phi(s,\cdot)$ is periodic in~$\overline{\Omega}$ for all $s\in\R$.
Lastly, since all functions~$u_n$ are increasing in time in
$\R\times\R^N$, the function~$u$ is such that~$u_t\ge 0$ in
$\R\times\overline{\Omega}$. From the previous observations and
parabolic regularity, there are then two periodic functions
$u^{\pm}$ defined in~$\overline{\Omega}$ such that $0\le u^-\le
u^+\le p_{\infty}$ in~$\overline{\Omega}$, $u(t,x)\to u^{\pm}(x)$ as
$t\to\pm\infty$ in $C^2_{loc}(\overline{\Omega})$ and~$u^{\pm}$ obey
$$\left\{\baa{rcll}
-\nabla\cdot(A(x)\nabla u^{\pm}) & = & F(x,u^{\pm}) & \hbox{in }\overline{\Omega},\vspace{3pt}\\
u^{\pm} & = & 0 & \hbox{on }\partial\Omega.\eaa\right.$$\par Let now
any index $i\in I_-$, that is~$\lambda_{1,\Omega_i,D}<0$ in the
sense of~(\ref{I-}). From the proof of Theorem~\ref{th2}, there is a
minimal periodic solution~$\widetilde{p}_i$ of~(\ref{statioi}).
Furthermore, in~$\overline{\Omega_i}$, there
holds~$\widetilde{p}_i=p\le p_{\infty}\le p_n$ for all $n\in\N$,
under the notation of Theorem~\ref{th3}. Therefore, one can always
shift in time the functions~$u_n$ so that, say,
$$\int_{C_0\cap\Omega_i}u_n(0,x)\,dx=\frac{1}{2}\int_{C_0\cap\Omega_i}\widetilde{p}_i(x)\,dx,$$
where we recall that $C_0=[0,L_1]\times\cdots\times[0,L_N]$. From
Lebesgue's dominated convergence theorem, the function $u$ satisfies
the same equality at the limit, whence
$$0\le\int_{C_0\cap\Omega_i}u^-(x)\,dx\le\frac{1}{2}\int_{C_0\cap\Omega_i}\widetilde{p}_i(x)\,dx
\le\int_{C_0\cap\Omega_i}u^+(x)\,dx$$ by monotonicity of~$u$ with
respect to~$t$. The minimality of~$\widetilde{p}_i$ and the strong
maximum principle imply that~$u^-=0$ in~$\overline{\Omega_i}$, while
$u^+>0$ in~$\Omega_i$, again from the strong maximum principle. If
we further assume that $F$ satisfies the KPP assumption~(\ref{Fkpp})
in~$\Omega$ (or just in~$\Omega_i$), then it follows from
Remark~\ref{rem31} that the solution of~(\ref{statioi}) is actually
unique, whence~$u^+=\widetilde{p}_i$ in $\overline{\Omega_i}$ in this case.\hfill\break

\noindent{\bf{Proof of part c) of Theorem~\ref{th4}.}} Firstly, it follows from~\cite{bh,w} that, for each $n\in\N$,
$$c^*_n(e)\ge\min_{\lambda>0}\frac{-k_{e,\lambda,n}}{\lambda}=\frac{-k_{e,\lambda_n,n}}{\lambda_n}$$
for some $\lambda_n>0$, where $k_{e,\lambda,n}$ denotes the
principal periodic eigenvalue of the operator \be\label{Lelambdan}
\mathcal{L}_{e,\lambda,n}:=-\nabla\cdot(A\nabla)+2\lambda\,Ae\cdot\nabla+\lambda\,\nabla\cdot(Ae)-\lambda^2\,
Ae\cdot e-\zeta_n\ \hbox{ in }\R^N. \ee Since, as above,
$k_{e,\lambda,n}\to k_{e,\lambda,D}$ as $n\to+\infty$
nondecreasingly for every $\lambda\in\R$, where $k_{e,\lambda,D}$ is
the principal periodic eigenvalue of the operator
$$\mathcal{L}_{e,\lambda,\Omega}:=-\nabla\cdot(A\nabla)+2\lambda\,Ae\cdot\nabla+\lambda\,
\nabla\cdot(Ae)-\lambda^2\,Ae\cdot e-\zeta\ \hbox{ in }\Omega$$
with zero Dirichlet boundary condition on~$\partial\Omega$, it follows that
$$c^*_n(e)\ge\frac{-k_{e,\lambda_n,n}}{\lambda_n}\ge\frac{-k_{e,\lambda_n,D}}{\lambda_n}
\ge\inf_{\lambda>0}\frac{-k_{e,\lambda,D}}{\lambda}$$ for all
$n\in\N$, whence \be\label{borneinf}
c^*(e)\ge\inf_{\lambda>0}\frac{-k_{e,\lambda,D}}{\lambda}. \ee
Furthermore, the maps $\lambda\mapsto-k_{e,\lambda,n}$ are all
convex and their derivatives at $\lambda=0$ are all equal to~$0$,
see~\cite{bh,bhr2}. In particular, for every $n\in\N$,
$-k_{e,\lambda,n}$ is nondecreasing with respect to~$\lambda\ge 0$
and~$-k_{e,\lambda,n}\ge-k_{e,0,n}=-\lambda_{1,n}$ for all
$\lambda\in\R$. By passing to the limit as~$n\to+\infty$ pointwise
in $\lambda$, one gets that the map $\lambda\mapsto-k_{e,\lambda,D}$
is convex in~$\R$, nondecreasing in $\R_+$, and there
holds~$-k_{e,\lambda,D}\ge-k_{e,0,D}=-\lambda_{1,D}>0$ for all
$\lambda\in\R$. Notice here that, if assumption~(\ref{discon}) is
made, then
$-k_{e,\lambda,D}\le\max_{\overline{\Omega}}\overline{\zeta}$ for
all $\lambda$, under the notation used in the proof of part~a).
Therefore, the infimum in~(\ref{borneinf}) is not reached in
general.\par Because of~(\ref{borneinf}), the
inequality~$-k_{e,\lambda,D}\ge-\lambda_{1,D}>0$ and the limit
$\lim_{n\to+\infty}k_{e,\lambda,n}=k_{e,\lambda,D}$ for
all~$\lambda$, it follows that, in order to show the positivity
of~$c^*(e)$, it is sufficient to prove that there exist $\Lambda>0$
and $\alpha>0$ such that \be\label{borneinfbis}
-k_{e,\lambda,n}\ge\alpha\,\lambda^2\ \hbox{ for all
}\lambda\ge\Lambda\hbox{ and for all }n\in\N. \ee Of course, from
the proof of part~a), this cannot be always true. However, assuming
from now on that~$A$ is constant, we shall now show
that~(\ref{borneinfbis}) holds under conditions~(\ref{slab})
or~(\ref{cylinder}). Assume first that there exist a unit vector
$e'\neq\pm e$ and two real numbers $a<b$ such that~(\ref{slab}) is
fulfilled, that is
$$\Omega\supset S_{e',a,b}:=\big\{x\in\R^N,\ a<x\cdot e'<b\big\}.$$
For any $\lambda>0$, let $\psi_{\lambda}$ be the function defined in $\overline{S_{e',a,b}}$ by
$$\psi_{\lambda}(x)=e^{\lambda'(x\cdot e')}\cos\Big(\frac{\pi}{b-a}\times\big(x\cdot e'-\frac{a+b}{2}\big)\Big),$$
where $\lambda'=\lambda\,(Ae\cdot e')/(Ae'\cdot e')$. The function
$\psi_{\lambda}$ is bounded and of class
$C^{\infty}\big(\overline{S_{e',a,b}}\big)$, it is positive in
$S_{e',a,b}$ and vanishes on $\partial S_{e',a,b}$. Furthermore,
since $\zeta_n=\zeta$ in
$\overline{\Omega}\supset\overline{S_{e',a,b}}$, it is
straightforward to check that
$$\mathcal{L}_{e,\lambda,n}\psi_{\lambda}=\Big(\frac{\pi^2(Ae'\cdot e')}{(b-a)^2}-\zeta(x)
-2\,\alpha\,\lambda^2\Big)\psi_{\lambda}\ \hbox{ in }\overline{S_{e',a,b}}$$
for all $n\in\N$, where $\alpha=(Ae\cdot e)/2-(Ae\cdot
e')^2/(2\,Ae'\cdot e')>0$ from Cauchy-Schwarz inequality, since the
unit vectors $e$ and $e'$ are not parallel. Since $\zeta$ is bounded
in $\overline{\Omega}$, it follows that there exists $\Lambda>0$ such that
$\mathcal{L}_{e,\lambda,n}\psi_{\lambda}\le-\alpha\,\lambda^2\,\psi_{\lambda}$
in $\overline{S_{e',a,b}}$ for all $\lambda\ge\Lambda$ and $n\in\N$.
This inequality yields~(\ref{borneinfbis}), as in the course of the
proof of Proposition~\ref{pro1}. We just sketch the proof here. Fix
any $\lambda\ge\Lambda$ and $n\in\N$ and let $\phi_n$ be a principal
periodic eigenfunction of the operator~$\mathcal{L}_{e,\lambda,n}$.
Namely, $\mathcal{L}_{e,\lambda,n}\phi_n=k_{e,\lambda,n}\phi_n$ and
$\phi_n$ is periodic and positive in~$\R^N$. Define
$$\epsilon^*=\sup\,\big\{\epsilon>0,\ \epsilon\psi_{\lambda}\le\phi_n\hbox{ in }\overline{S_{e',a,b}}\big\}.$$
Owing to the definition of~$\psi_{\lambda}$ and the uniform
positivity of $\phi_n$, the quantity~$\epsilon^*$ is a posi\-tive
real number. Furthermore, $\epsilon^*\psi_{\lambda}\le\phi_n$ in
$\overline{S_{e',a,b}}$ and there is a sequence~$(x_m)_{m\in\N}$ of points in~$S_{e',a,b}$ such that
$\liminf_{m\to+\infty}d(x_m,\partial S_{e',a,b})\!>\!0$,
$\lim_{m\to+\infty}(\epsilon^*\psi_{\lambda}(x_m)\!-\!\phi_n(x_m))\!=\!0$
and~$\liminf_{m\to+\infty}\mathcal{L}_{e,\lambda,n}(\epsilon^*\psi_{\lambda}-\phi_n)(x_m)\ge0$.
Since there holds $\mathcal{L}_{e,\lambda,n}\phi_n(x_m)=k_{e,\lambda,n}\phi_n(x_m)$ and
$\mathcal{L}_{e,\lambda,n}\psi_{\lambda}(x_m)\le-\alpha\lambda^2\psi_{\lambda}(x_m)$
for every $m\in\N$, one concludes that
$k_{e,\lambda,n}\le-\alpha\lambda^2$, that is~(\ref{borneinfbis}).
This yields the desired inequality~$c^*(e)>0$, as already
emphasized.\par Assume now that there exist a unit vector $e'$, a
point $x_0\in\R^N$ and a positive real number~$r$ such that $e'$ is
an eigenvector of $A$ with $e'\cdot e\neq 0$ and~(\ref{cylinder})
holds. Let $\beta>0$ be such that~$Ae'=\beta e'$. Since the matrix
$A$ is symmetric, there is an orthonormal family of eigenvectors
$e'_1,\ldots,e'_{N-1}$ of~$A$ in~$\R^N$ such that~$e'_i\cdot e'=0$
for all $1\le i\le N-1$. Even if it means decreasing $r>0$
in~(\ref{cylinder}), one can assume without loss of generality that
$$\Omega\supset C_{e',r}:=\big\{x\in\R^N,\ |(x-x_0)\cdot e'_i|<r\hbox{ for all }1\le i\le N-1\big\}.$$
For any $\lambda>0$, let $\psi_{\lambda}$ be the function defined in $\overline{C_{e',r}}$ by
$$\psi_{\lambda}(x)=\prod_{1\le i\le N-1}e^{\lambda'_i(x\cdot e'_i)}\cos\Big(\frac{\pi(x-x_0)\cdot e'_i)}{2\,r}\Big),$$
where $\lambda'_i=\lambda\,(Ae\cdot e'_i)/(Ae'_i\cdot e'_i)$. The
function $\psi_{\lambda}$ is bounded and of class
$C^{\infty}\big(\overline{C_{e',r}}\big)$, it is positive in
$C_{e',r}$ and vanishes on $\partial C_{e',r}$. Furthermore, since
$\zeta_n=\zeta$ in $\overline{\Omega}\supset\overline{C_{e',r}}$, it
is straightforward to check that
$$\mathcal{L}_{e,\lambda,n}\psi_{\lambda}=\Big[\big(\sum_{1\le i\le N-1}\!\!\!\frac{\pi^2(Ae'_i\cdot e'_i)}{4r^2}\big)
-\zeta(x)-2\,\alpha\,\lambda^2\Big]\psi_{\lambda}\ \hbox{ in
}\overline{C_{e',r}}$$ for all $n\in\N$, where $\alpha=\beta(e\cdot
e')^2/2>0$ since $\beta>0$ and $e\cdot e'\neq 0$ by assumption.
Thus, one concludes as above that there is $\Lambda>0$ such that
$\mathcal{L}_{e,\lambda,n}\psi_{\lambda}\le-\alpha\,\lambda^2\psi_{\lambda}$
in $\overline{C_{e',r}}$ for all $\lambda\ge\Lambda$ and $n\in\N$.
This yields~(\ref{borneinfbis}) and finally $c^*(e)>0$.
The proof of Theorem~\ref{th4} is thereby complete.\hfill$\Box$

\begin{rem}{\rm In the case when the functions $f_n$ fulfill the KPP assumption~$(\ref{fkppn})$,
then~$c^*(e)$ is given by an explicit variational formula. Namely,
under assumption~$(\ref{fkppn})$ for the functions~$f_n$, it follows
from the proof of part~{\rm{a)}} of Theorem~$\ref{th4}$ with the
choices $\overline{f}_n=f_n$ and~$\overline{F}=F$
that~$c^*(e)\le\inf_{\lambda>0}-k_{e,\lambda,D}/\lambda$, because
of~$(\ref{c*e})$. On the other hand, the reverse
inequality~$(\ref{borneinf})$ always holds, from the proof of
part~c) of Theorem~$\ref{th4}$. As a conclusion, the
assumption~$(\ref{fkppn})$ for the functions $f_n$ yields
$$c^*(e)=\inf_{\lambda>0}\frac{-k_{e,\lambda,D}}{\lambda}.$$}
\end{rem}


\footnotesize{

}


\begin{thebibliography}{ABC}
\bibitem{aw} D.G. Aronson and H.F. Weinberger, Multidimensional nonlinear diffusions arising in population genetics,
{\it Adv. Math.} {\bf 30} (1978), 33-76.
\bibitem{b} H. Berestycki, Le nombre de solutions de certains probl\`emes semi-lin\'eaires elliptiques,
{\it J.~Funct. Anal.} {\bf 40} (1981), 1-29.
\bibitem{bh} H. Berestycki and F. Hamel, Front propagation in periodic excitable media,
{\it Comm. Pure Appl. Math.} {\bf 55} (2002), 949-1032.
\bibitem{bhna} H. Berestycki, F. Hamel and G. Nadin, Asymptotic spreading in heterogeneous diffusive media,
{\it J.~Funct. Anal.} {\bf 255} (2008), 2146-2189.
\bibitem{bhr1} H. Berestycki, F. Hamel and L. Roques, Analysis of the periodically fragmented environment model: I~--~Species persistence,
{\it J.~Math. Biol.} {\bf 51} (2005), 75-113.
\bibitem{bhr2} H. Berestycki, F. Hamel and L. Roques, Analysis of the periodically fragmented environment model:
II~--~Biological invasions and pulsating travelling fronts, {\it J.~Math. Pures Appl.} {\bf 84} (2005), 1101-1146.
\bibitem{bhr3} H.~Berestycki, F. Hamel and L.~Rossi, Liouville type results for semilinear elliptic equations in unbounded domains,
{\it Ann. Mat. Pura Appl.} {\bf 186} (2007), 469-507.
\bibitem{bn} H. Berestycki and L. Nirenberg, Travelling fronts in cylinders, {\it Ann. Inst. H.~Poincar\'e, Analyse Non Lin\'eaire} {\bf 9} (1992), 497-572.
\bibitem{bnv} H. Berestycki, L. Nirenberg and S.R.S. Varadhan, The principal eigenvalue and maximum principle for second order elliptic operators in general domains,
{\it Comm. Pure Appl. Math.} {\bf 47} (1994), 47-92.
\bibitem{cc1} R.S. Cantrell and C. Cosner, The effects of spatial heterogeneity in population dynamics, {\it J. Math. Biol.} {\bf 29} (1991), 315-338.
\bibitem{cc2} R.S. Cantrell and C. Cosner, On the effects of spatial heterogeneity on the persistence of interacting species,
{\it J.~Math. Biol.} {\bf 37} (1998), 103-145.
\bibitem{dm} Y. Du and H. Matano, Convergence and sharp thresholds for propagation in nonlinear diffusion problems, {\it J.~Europ. Math. Soc.} {\bf 12} (2010), 279-312.
\bibitem{e1} M. El Smaily, Pulsating travelling fronts: Asymptotics and homogenization regimes, {\it Europ.~J. Appl. Math.} {\bf 19} (2008), 393-434.
\bibitem{f} R.A. Fisher, The advance of advantageous genes, {\it Ann. Eugenics} {\bf 7} (1937), 335-369.
\bibitem{fg} M. Freidlin and J. G\"artner, On the propagation of concentration waves in periodic and random media,
{\it Sov. Math. Dokl.} {\bf 20} (1979), 1282-1286.
\bibitem{h8} F. Hamel, Qualitative properties of monostable pulsating fronts: exponential decay and monotoni\-city,
{\it J.~Math. Pures Appl.} {\bf 89} (2008), 355-399.
\bibitem{hr} F. Hamel and L.~Roques, Uniqueness and stability properties of monostable pulsating fronts,
{\it J.~Europ. Math. Soc.} {\bf 13} (2011), 345-390.
\bibitem{h} S. Heinze, Large convection limits for KPP fronts, {\it preprint}.
\bibitem{kks} N. Kinezaki, K. Kawasaki and N. Shigesada, Spatial dynamics of invasion in sinusoidally varying environments,
{\it Popul. Ecol.} {\bf 48} (2006), 263-270.
\bibitem{kpp} A.N. Kolmogorov, I.G. Petrovsky and N.S. Piskunov, \'Etude de l'\'equation de la diffusion avec croissance de la quantit\'e de mati\`ere et son
application \`a un probl\`eme biologique, {\it Bull. Univ. Etat Moscou, S\'er. Intern.~A} {\bf 1} (1937), 1-26.
\bibitem{llm} X. Liang, X. Lin and H. Matano, A variational problem associated with the minimal speed of travelling waves for spatially periodic
reaction-diffusion equations, {\it Trans. Amer. Math. Soc.} {\bf 362} (2010), 5605-5633.
\bibitem{lz1} X. Liang and X.Q. Zhao, Asymptotic speeds of spread and traveling waves for monostable semiflows with applications,
{\it Comm. Pure Appl. Math.} {\bf 60} (2007), 1-40.
\bibitem{lz2} X. Liang and X.Q. Zhao, Spreading speeds and traveling waves for abstract monostable evolution systems,
{\it J.~Funct. Anal.} {\bf 259} (2010), 857-903.
\bibitem{mn} C.B. Muratov and M. Novaga, Front propagation in infinite cylinders. I.~A variational approach, {\it Comm. Math. Sci.} {\bf 6} (2008), 799-826.
\bibitem{mu} J.D. Murray, {\it Mathematical Biology}, Springer, 2003.
\bibitem{n1} G. Nadin, Travelling fronts in space-time periodic media, {\it J.~Math. Pures Appl.} {\bf 92} (2009), 232-262.
\bibitem{n4} G. Nadin, Some dependence results between the spreading speed and the coefficients of the space-time periodic Fisher-KPP equation,
{\it Europ. J.~Appl. Math.} (2011), to appear.
\bibitem{nrx} J. Nolen, M. Rudd and J. Xin, Existence of KPP fronts in spatially-temporally periodic advection and variational principle for propagation speeds,
{\it Dyn. Part. Diff. Equations} {\bf 2} (2005), 1-24.
\bibitem{nx} J. Nolen and J. Xin, Existence of KPP type fronts in space-time periodic shear flows and a study of minimal speeds based on variational principle,
{\it Disc. Cont. Dyn. Systems~A} {\bf 13} (2005), 1217-1234.
\bibitem{p} P. Polacik, Threshold solutions and sharp transitions for nonautonomous parabolic equations on $\R^N$, {\it Arch. Ration. Mech. Anal.} {\bf 199} (2011), 69-97.
\bibitem{r} J.-M. Roquejoffre, Eventual monotonicity and convergence to travelling fronts for the solutions of parabolic equations in cylinders,
{\it Ann. Inst. H.~Poincar\'e, Analyse Non Lin\'eaire} {\bf 14} (1997), 499-552.
\bibitem{rz} L. Ryzhik and A. Zlato{\v{s}}, KPP pulsating front speed-up by flows, {\it Comm. Math. Sci.} {\bf 5} (2007), 575-593.
\bibitem{sk}  N. Shigesada and K. Kawasaki, {\it Biological Invasions: Theory and Practice}, Oxford Series in Ecology and Evolution, Oxford: Oxford Univ. Press, 1997.
\bibitem{v2} J.M. Vega, Travelling waves fronts of reaction-diffusion equations in cylindrical domains, {\it Comm. Part. Diff. Equations} {\bf 18} (1993), 505-531.
\bibitem{w} H.F. Weinberger, On spreading speeds and traveling waves for growth and migration in periodic habitat, {\it J. Math. Biol.} {\bf 45} (2002), 511-548.
\bibitem{x3} X. Xin, Existence of planar flame fronts in convective-diffusive periodic media, {\it Arch. Ration. Mech. Anal.} {\bf 121} (1992), 205-233.
\bibitem{x6} J.X. Xin, Analysis and modeling of front propagation in heterogeneous media, {\it SIAM Review} {\bf 42} (2000), 161-230.
\bibitem{z1} A. Zlato{\v{s}}, Sharp transition between extinction and propagation of reaction, {\it J.~Amer. Math. Soc.} {\bf 19} (2006), 251-263.
\bibitem{z2} A. Zlato{\v{s}}, Sharp asymptotics for KPP pulsating front speed-up and diffusion enhancement by flows,
{\it Arch. Ration. Mech. Anal.} {\bf 195} (2010), 441-453.
\end{thebibliography}
\end{document}